\documentclass[12pt]{amsart}
\usepackage{cases}
\usepackage{txfonts}
\usepackage{latexsym}
\usepackage{amsthm}
\usepackage{mathrsfs}
\usepackage{amssymb, amsmath}
\usepackage{cite}

\newtheorem{theorem}{Theorem}

\newtheorem{lemma}[theorem]{Lemma}

\numberwithin{equation}{section}
\numberwithin{theorem}{section}

\newtheorem{proposition}[theorem]{Proposition}

\def\qed{\ifhmode\textqed\fi
   \ifmmode\ifinner\quad\qedsymbol\else\dispqed\fi\fi}
\def\textqed{\unskip\nobreak\penalty50
    \hskip2em\hbox{}\nobreak\hfil\qedsymbol
    \parfillskip=0pt \finalhyphendemerits=0}
\def\dispqed{\rlap{\qquad\qedsymbol}}

\begin{document}
\title{T\MakeLowercase{he twisted mean square and critical zeros of }D\MakeLowercase{irichlet} $L$-\MakeLowercase{functions}}
\author{Xiaosheng Wu}
\date{}
\address {School of Mathematics, Hefei University of Technology, Hefei 230009,
P. R. China.}
\email {xswu@amss.ac.cn}
\thanks{This work is supported by the National Natural Science Foundation of China (Grant No. 11871187) and the Fundamental Research Funds for the Central Universities.}

\begin{abstract} In this work, we obtain an asymptotic formula for the twisted mean square of a Dirichlet $L$-function with a longer mollifier, whose coefficients are also more general than before.  As an application we obtain that, for every Dirichlet $L$-function, more than 41.72\% of zeros are on the critical line and more than 40.74\% of zeros are simple and on the critical line. These proportions also improve previous results which were proved only for the Riemann zeta-function.
\end{abstract}
\maketitle

\section{Introduction}
Let $\chi$ be a Dirichlet character with $q$ its modulus and $L(s,\chi)$ be its associated Dirichlet $L$-function. When $s=\sigma+it$ with $\sigma>1$, we define $L(s,\chi)$  by
\begin{align}
   L(s,\chi)=\sum_{n\ge1}\chi(n)n^{-s}.
\end{align}
We are interested in an asymptotic formula for
\begin{align}
\label{DefI}
   I(\chi)=\int_T^{2T}\bigg|L\bigg(\frac12+it,\chi\bigg)\bigg|^2 \bigg|B\bigg(\frac12+it,\chi\bigg)\bigg|^2dt,
\end{align}
where  $\chi$ is a primitive Dirichlet character and $B(s,\chi)$ is a Dirichlet polynomial
\begin{align}
   B(s,\chi)=\sum_{n\le y}\frac{\chi(n)a(n)}{n^s}\ \ \ \ \text{with} \ \ \ \ a_n\ll n^\epsilon, \ \ \ y=T^\theta,\ \ \ \text{and}\ \ \ \theta<1.
\end{align}
We restrict $\chi$ to a primitive Dirichlet character since most properties of Dirichlet $L$-functions to non-primitive characters can be deduced directly from ones to corresponding primitive characters.

Asymptotic formulae for $I$ have been widely applied in studying $L$-functions, especially in the distribution of values of $L$-functions, the location of their critical zeros and upper and lower bounds for the size of $L$-functions. See for example, \cite{BG, Con-More, CIS-Crit, CGG,Rad,Sou}.

The value of $\theta$ is crucially important since it limits the best result we can obtain in most cases. For example, a larger $\theta$ means a larger lower bound for the proportion of critical zeros, and  $\theta=1-\epsilon$ means the Lindel\"{o}f hypothesis. Moreover, Bettin and Gonek \cite{BG} have proved that $\theta=\infty$ implies the Riemann hypothesis while it is normally conjectured $\theta<1$.

An asymptotic formula for $I$ was firstly obtained for the Riemann zeta-function.  In 1985, Balasubramanian, Conrey and Heath-Brown in \cite{BCH} gave an explicit formula for the Riemann zeta-function that
\begin{align}\label{+1.6}
I(1)=T\sum_{h,k\le y}\frac{a(h)\overline{a(k)}}{[h,k]}\bigg(\log\frac{T(h,k)^2}{2\pi hk}+2\gamma-1+2\log2\bigg)+O(T^{1-\epsilon_\theta}),
\end{align}
where $\epsilon_\theta$ was a constant decided by $\theta$, and $(h,k)$, $[h,k]$ denoted the gcd and the lcm of $h$, $k$, respectively. In general, they proved that $\epsilon_\theta>0$ with $\theta<\frac12$. Further, for a special form
\begin{align}\label{++1}
a(m)=\mu(m)\mathcal{F}(m)\ \ \ \text{with}\ \ \ \mathcal{F}\in\mathbb{F}=\bigg\{\mathcal{F}:\mathcal{F}(x)\ll_\epsilon x^\epsilon,~ \mathcal{F}'(x)\ll\frac1x\bigg\},
\end{align}
they obtained $\epsilon_\theta>0$ with $\theta<\frac9{17}$. With this larger value of $\theta$ they improved the proportion of zeros of the Riemann zeta-function on the critical line from at least 36.58\% to 38\%.\\

\noindent
{\bf Definition} [Separability].
Let $\mathbb{S}$ be a set of arithmetical functions. We say that $\mathcal{F}$ is separable or has a property of separability in $\mathbb{S}$ if $\mathcal{F}\in \mathbb{S}$ and $\mathcal{F}(mn)$ can be expressed as a finite sum of $\mathcal{F}_i(m)\mathcal{F}_j(n)$ with $\mathcal{F}_i,\mathcal{F}_j\in \mathbb{S}$.\\

Let $\chi$ be a Dirichlet character with $q$ its modulus, we denote
\begin{align}
\mathcal{L}_\chi:=\log\frac{qT}{2\pi}\ \ \ \ \text{and}\ \ \ \ \mathcal{L}_1:=\log\frac{T}{2\pi}.
\end{align}
Let $Q(x)$ be a polynomial,
\begin{align}\label{x1.7}
V(s,\chi)=Q\bigg(-\frac1{\mathcal{L}_\chi}\frac{d}{ds}\bigg)L(s,\chi)
\end{align}
and
\begin{align}
\label{x1.6}
I_R(Q,\chi)=\int_T^{2T}\bigg|V\bigg(\frac12+it+\frac{R}{\mathcal{L}_\chi},\chi\bigg)\bigg|^2 \bigg|B\bigg(\frac12+it,\chi\bigg)\bigg|^2dt
\end{align}
for any given real number $R$.
In 1989, Conrey \cite{Con-More} proved an asymptotic formula for $I_R(Q,1)$. Conrey worked with a coefficient similar as \eqref{++1} but with $\mathcal{F}$ separable in $\mathbb{F}$. By employing Weil's bound and a work of Deshouillers and Iwaniec \cite{DIK,DIP} on Kloosterman sums to control error terms, Conrey extended $\theta$ to $\frac4{7}$ and proved that the Riemann zeta-function has more than 40.88\% of zeros on the critical line.

Let $\Delta=T^{1-\eta}$ with a small constant $\eta>0$, Conrey \cite{Con-More} actually deduced the asymptotic formula for $I_R(Q,1)$ from an asymptotic formula of
\begin{align}
\int_{-\infty}^{\infty}e^{-(t-T)^2\Delta^{-2}}\zeta\bigg(\frac12+it+\alpha\bigg)\zeta\bigg(\frac12-it+\beta\bigg) \bigg|B\bigg(\frac12+it,1\bigg)\bigg|^2dt
\end{align}
for any $\alpha,~\beta\ll \mathcal{L}^{-1}_1$. The factor $e^{-(t-T)^2\Delta^{-2}}$ here is actually a smooth function to smoothen the integral. In 2010, Young \cite{You} deduced an asymptotic formula for $I_R(Q,1)$ by considering an asymptotic formula of
\begin{align}
\label{x1.9}
\int_{-\infty}^{\infty}\Phi\bigg(\frac tT\bigg)\zeta\bigg(\frac12+it+\alpha\bigg)\zeta\bigg(\frac12-it+\beta\bigg) \bigg|B\bigg(\frac12+it,1\bigg)\bigg|^2dt,
\end{align}
where $\Phi(x)$ was a smooth function that should satisfy some properties to smoothen the integral. With $\theta<\frac12$ Young deduced an asymptotic formula for \eqref{x1.9} in a different way, which involved an approximate functional equation of $\zeta(\frac12+\alpha+it)\zeta(\frac12+\beta-it)$ and had a shorter proof.

With $\alpha=\beta=0$ and any coefficient $a(m)\ll m^\epsilon$ in \eqref{x1.9}, Bettin, Chandee and Radziwi{\l}{\l} \cite {BCR} broke the $\theta<\frac12$ barrier. Specifically, an asymptotic formula with $\theta<\frac{17}{33}$ is proved in \cite{BCR}. Their work is an extension of Young's approach, which begins with an approximate function equation of $|\zeta(1/2+it)|^2$ pointed out by Li and Radziwi{\l}{\l} in \cite{LR}. Then an estimate for trilinear forms of Kloosterman fractions in \cite{BC} is vital to their estimation on error terms. They also proved some interesting applications, especially an upper bound for $2k$-th moment of the Riemann zeta-function with $k=1+1/n$.

Recently, Pratt and Robles \cite{PR} proved an asymptotic formula for \eqref{x1.9} by extending the way in \cite{BCR} and \cite{You}. In more specific terms, an asymptotic formula was obtained when $\theta<\frac{17}{33}$ for general coefficients $a(n)\ll n^\epsilon$ and $\theta<\frac47$ for a special coefficient as in \cite{Con-More}. In addition, they also extended $\theta$ to $\frac6{11}$ for the special coefficient $a(n)=\mu^2(n)(\mu*\Lambda^{*k})(n)P(\frac{\log y/n}{\log y})$ in Feng's mollifier, and thus improving the proportion of zeros of the Riemann zeta-function on the critical line to at least 41.491\%.

In contrast to the Riemann zeta-function, results on general Dirichlet $L$-functions are weaker. An asymptotic formula of $I_R(Q,\chi)$ was proved only when $\theta<\frac12$ and $q=o(\log T)$ by Bauer \cite{Bau} in 2000. Bauer found that a Kloosterman sum in error terms with a Dirichlet character was too complicated to estimate. Thus, he selected a different way to avoid this Klooserman sum. He actually followed a way used in Conrey and Ghosh \cite{CG}, which can avoid the Kloosterman sum by some large sieve inequalities but only worked with $\theta<\frac12$.

It does not seem that we can avoid Kloosterman sums when extend $\theta$ to the right-hand side of $\frac12$. However, the estimation on a Kloosterman sum with a Dirichlet character seems to be very complicated. A possible cause is the discreteness of $\chi$, which makes both coefficients of $L(s,\chi)$ and $B(s,\chi)$ non-differentiable, while special forms of $a(m)$ in \cite{BCH}, \cite{Con-More} and \cite{PR} rely heavily on the differentiability of $\mathcal{F}$. In this work we try to extend $\theta$ to the right-hand side of $\frac12$ by the way of \cite{BCH,Con-More}. We note that the way developed in \cite{You}, also used in \cite{Bau} and \cite{PR}, makes the calculation for  the Rieman zeta-function `cleaner', but it does not seem to be so for $L(s,\chi)$. Because of the existence of a Dirichlet character, the calculation seems to be much more complicated in this way, especially when one tries to distinguish character information from off-diagonal terms. For a Dirichlet $L$-function, we will extend $\theta$ to the right-hand side of $\frac12$ uniformly in $q$ with $\log q=o(\log T)$ for different forms of coefficients. Our estimation on error terms is also based on Weil's bound, Deshouillers and Iwaniec's estimate on Kloosterman sums and the estimate of trilinear forms of Kloosterman fractions obtained in \cite{BC}. However, before using these estimates, we should make some technical preparation carefully to strip its entanglement with the Dirichlet character first. In addition, we obtain a more general form of $a(m)$ that can be applied to $I_R(Q,\chi)$ with $\theta<\frac47$. With a coefficient of this form, we obtain a larger lower bound for the proportion of zeros on the critical line.

In this work, we give an asymptotic formula for $I_R(Q,\chi)$, and our results are specified in the following theorem.

\begin{theorem}\label{thm1}
Let $\chi$ be a primitive Dirichlet character (mod $q$) with $\log q=o(\log T)$ and $\alpha=a/\mathcal{L}_\chi,~\beta=b/\mathcal{L}_\chi$ with $a,~b\in\mathbb{C}$ and $a,~b\ll1$. Let $I_R(Q,\chi)$ be defined as in \eqref{x1.6}. Suppose that $a(m)\ll_\epsilon m^{\epsilon}$ for any $\epsilon>0$ and $y=T^\theta$, then we have
\begin{align}\label{+1.10}
I_R(Q,\chi)&=TQ\bigg(\frac{-d}{da}\bigg)\overline{Q}\bigg(\frac{-d}{db}\bigg) \Bigg\{\sum_{h,k\le y}\frac{\chi_0(hk)(h,k)^{\alpha+\beta}a(h)\overline{a(k)}}{h^{1+\beta} k^{1+\alpha}}\\
&\times\Bigg(\frac{2^{1+\alpha+\beta}-1}{1+\alpha+\beta}\bigg(\frac {2\pi hk}{qT(h,k)^2}\bigg)^{\alpha+\beta}L(1-\alpha-\beta,\chi_0) +L(1+\alpha+\beta,\chi_0)\Bigg)\Bigg\}\Bigg|_{a=b=-R}+O(T^{1-\epsilon_\theta}),\notag
\end{align}
and, in the particular case $\alpha=\beta=0$,
\begin{align}\label{+1.12}
I(\chi)=T\frac{\phi(q)}{q}\sum_{h,k\le y}\frac{a(h)\overline{a(k)}}{[h,k]}\chi_0(hk)\bigg(\log\frac{Tq(h,k)^2}{2\pi HK}+2\gamma-1+c_q+2\log2)\bigg)+O(T^{1-\epsilon_\theta})
\end{align}
with $c_q=\sum_{p\mid q}(\log p)/(p-1)$ and $\gamma$ is the Euler's constant. Here $\epsilon_\theta$ is a constant depending on $\theta$ as follows:
\begin{description}
  \item[(A)] In general, we have $\epsilon_\theta>0$ for any given $\theta<\frac{17}{33}$;
  \item[(B)] We have $\epsilon_\theta>0$ for any given $\theta<\frac{4}{7}$ when $a(n)$ has a special form
\begin{align}
a(n)=\mu(n)(\mathcal{F}_0+\mathcal{F}_1\cdot(\mathcal{F}_2*\mathcal{F}_3))(n)\notag
\end{align}
with $\mathcal{F}_i$ separable in $\mathbb{F}=\big\{\mathcal{F}:\mathcal{F}(x)\ll_\epsilon x^\epsilon,~ \mathcal{F}'(x)\ll\frac1x\big\}$ for $0\le i\le 3$. In addition, it also holds when one of $\mathcal{F}_2$ and $\mathcal{F}_3$ is separable in $\big\{\mathcal{F}:\mathcal{F}(x)\ll_\epsilon x^\epsilon,~\mathcal{F}(x)=0\ \ \text{for}\ \ x> y^{\frac34}\big\}$ and other $\mathcal{F}_i$ are separable in $\mathbb{F}$.
\end{description}
\end{theorem}

\noindent
{\bf Remark.} We give some remarks as follows:
\begin{itemize}
 \item Theorem \ref{thm1} only gives the primitive character case, however, if wanted, one can deduce similar results for non-primitive characters since the restriction of primitive characters is only used to simplify coefficients in our proof.
  \item Similar results on upper bounds for $2k$-th moment of Dirichlet $L$-functions  in $t$-aspect with $k=1+1/n$ can be obtained in the same way as \cite{BCR} by (A), these upper bounds should be uniform in $q$ with $\log q=o(\log T)$.
\end{itemize}

When we apply (B) of Theorem \ref{thm1}, we obtain that every Dirichlet $L$-function has more than 41.72\% of zeros on the critical line. We present it in Theorem \ref{thm2} and give its proof in Section \ref{sec5}.

Let $N(T,\chi)$ denote the number of zeros of $L(s,\chi)$ with $0<\sigma<1$ and $|t|\le T$. Also let $N_c(T,\chi)$ and $N^*_c(T,\chi)$ denote the number of zeros and simple zeros of $L(1/2+it, \chi)$ with $|t|\le T$ respectively. Then we define $\kappa(\chi)$ and $\kappa^*(\chi)$ by
\begin{align}
\kappa(\chi)=\frac{N_c(T,\chi)}{N(T,\chi)},\ \ \ \ \kappa^*(\chi)=\frac{N^*_c(T,\chi)}{N(T,\chi)}.
\end{align}

There is a long history of studying on the proportion of zeros lying on the critical line for the Riemann zeta-function, and one may see \cite{Con-More,BCY,Feng,Lev-More} for example. By the approach of Levinson \cite{Lev-More}, Conrey \cite{Con-More} and the observation of Heath-Brown \cite{Hea}, it is known that \cite{PR,Bui}
\begin{align}
\kappa(1)>.41491\ \ \ \ \ and \ \ \ \ \ \kappa^*(1)>.40589
\end{align}
for sufficiently large $T$.

For a general Dirichlet $L$-function, it is proved by Bauer \cite{Bau} in 2000 that
\begin{align}
\kappa(\chi)>.365815\ \ \ \ \ and \ \ \ \ \ \kappa^*(\chi)>.356269
\end{align}
for sufficiently large $T$ with  $q=o(\log T)$ .

\begin{theorem}\label{thm2}
We have, for any Dirichlet character $\chi$,
\begin{align}
\kappa(\chi)>.4172\ \ \ \ \ and \ \ \ \ \ \kappa^*(\chi)>.4074
\end{align}
for sufficiently large $T$ with $\log q=o(\log T)$.
\end{theorem}

\noindent
{\bf Remark.}
We do not need to restrict this theorem to primitive characters since Dirichlet $L$ functions to non-primitive characters share the same non-trivial zeros as ones to corresponding primitive characters.\\

It should also be noted that some much further results have been obtained for the family of Dirichlet $L$-functions. By averaging over all primitive characters and all $q\le Q$ with $Q$ restricted to be sufficiently large in  terms of $T$, Conrey, Iwaniec and Soundararajan \cite{CIS-Crit} proved that at least 58.65\% of zeros  of the family of Dirichlet $L$-functions are simple and on the cirtical line, and the percentage only in the simple zeros case has been improved to 60.261\% by Wu \cite{Wu}. These results rely on Conrey, Iwaniec and Soundararajan's Asymptotic Large Sieve work \cite{CIS-Asym}.

To prove Theorem \ref{thm2}, we use a mollifier
\begin{align}\label{xm1.18}
\psi(s)=\sum_{n\le y}\frac{\chi(n)a(n)}{n^{s+R/\mathcal{L}_\chi}}
\end{align}
with its coefficient
\begin{align}\label{xa1.18}
a(n)=\mu(n)\bigg(P_1\bigg(\frac{\log y/n}{\log y}\bigg)+P_2\bigg(\frac{\log y/n}{\log y}\bigg)\sum_{p\mid n,~p\le y^{3/4}}P\bigg(\frac{\log p}{\log y}\bigg)\bigg).
\end{align}
Here $P_i$ and $P$ are real polynomials that satisfy some minor conditions. This coefficient can be seen as two parts, the first one is due to Conrey's mollifier and the second one is motivated by Feng's mollifier. Different to Feng's mollifier, an interesting mollifier of two pieces
\begin{align}
\psi(s)=\sum_{n\le y_1}\frac{P_1\Big(\frac{\log y_1/n}{\log y_1}\Big)}{n^{s+R/\mathcal{L}_1}}+\chi(s-R/\mathcal{L}_1)\sum_{hk\le y_2}\frac{(\mu*\mu)(h)P_2\Big(\frac{\log y_2/hk}{\log y_2}\Big)}{h^{s+R/\mathcal{L}_1}k^{1-s-R/\mathcal{L}_1}}
\end{align}
with $y_1\le T^{\frac47}$ and $y_2\le T^{\frac12}$ is also introduced by \cite{BCY}, and some extensions have also been obtained in \cite{Bui,KRZ,Son, RRZ}.

 Let us see \eqref{xa1.18} in detail. In Feng's mollifier, an additional part, which was deduced from
\begin{align}\label{+1.18}
\sum_{2\le k\le K}\frac1{\log^k y}(\mu*\underbrace{\Lambda*\Lambda*\cdots*\Lambda}_k)(n)P_k\bigg(\frac{\log y/n}{\log y}\bigg),
\end{align}
was injected into the coefficient. If we ignore non-squarefree $n$ terms, Feng's mollifier can be thought as a simplification of `continuous' truncation of the Dirichlet series
\begin{align}
\label{x1.22}
\frac1{\zeta(s)+\frac{\zeta'(s)}{\mathcal{L}_1}}=\frac1{\zeta(s)}-\frac{\zeta'(s)} {\mathcal{L}_1\zeta^2(s)}+\frac{\zeta'^2(s)}{\mathcal{L}^2_1\zeta^3(s)}- \frac{\zeta'^3(s)}{\mathcal{L}^3_1\zeta^4(s)}+\cdots.
\end{align}
Let us distinguish the major contributor in Feng's coefficient. We only consider squarefree $n$ as Feng did, then \eqref{+1.18} is equal to
\begin{align}
(-1)^k\mu(n)\frac1{\log^ky}(1*\underbrace{\Lambda*\Lambda*\cdots*\Lambda}_k)(n)P_k\bigg(\frac{\log y/n}{\log y}\bigg).
\end{align}
It is easy to see that
\begin{align}
\mu(n)\frac1{\log^ky}(1*\underbrace{\Lambda*\Lambda*\cdots*\Lambda}_k)(n)&=\mu(n)\sum_{p_1\cdots p_k\mid n}\frac{\log p_1\cdots\log p_k}{\log^k y}\notag\\
&=\mu(n)\sum_{p_1\cdots p_{k-1}\mid n}\frac{\log p_1\cdots\log p_{k-1}}{\log^{k} y}\sum_{p_k\mid n/p_1\cdots p_{k-1}}\log p_k\notag\\
&=\mu(n)\sum_{p_1\cdots p_{k-1}\mid n}\frac{\log p_1\cdots\log p_{k-1}}{\log^{k} y}\log\Bigg(\frac{\log n}{\log (p_1\cdots p_{k-1})}\Bigg)
\end{align}
by the simple formula $\sum_{p\mid n}\log p=\log n$ for squarefree $n$. This splits $\mu(n)\frac1{\log^ky}(1*\underbrace{\Lambda*\Lambda*\cdots*\Lambda}_k)(n)$ into following two terms
\begin{align}
\mu(n)\frac{\log n}{\log y}\sum_{p_1\cdots p_{k-1}\mid n}\frac{\log p_1\cdots\log p_{k-1}}{\log^{k} y}-(k-1)\mu(n)\sum_{p_1\cdots p_{k-1}\mid n}\frac{\log^2 p_1\cdots\log p_{k-1}}{\log^{k-1} y}.
\end{align}
Then we can use the formula $\sum_{p\mid n}\log p=\log n$ again to eliminate $p_{k-1}$ in both terms. By repeating this action we may eliminate all factors $\log p_i$ with degree one, and ultimately obtain

\begin{align}\label{x1.30}
\mu(n)\frac1{\log^ky}(1*&\underbrace{\Lambda*\Lambda*\cdots*\Lambda}_k)(n)\notag\\
& \ \ \ \ \ \ \ \ \ \ \ \ \ \ \ =\mu(n)\sum_{j=0}^{k/2}\sum_{{k_0+k_1+\cdots+k_j=k}\atop{k_1,\cdots,k_j\ge2}} a_{k_1,\cdots,k_j}\Bigg(\frac{\log n}{\log y}\Bigg)^{k_0} \sum_{p_1\cdots p_j\mid n}\frac{\log^{k_1} p_1\cdots\log^{k_j} p_j}{\log^{k-k_0}y},
\end{align}
where $a_{k_1,\cdots,k_j}$ are some combinatorial constants, to see \cite{PRZZ} for their exact values.
As Feng \cite{Feng} did, we only consider $k\le3$ here for convenience. It is obvious that $j$ only takes values $0$ and $1$ in \eqref{x1.30} with $k\le3$. In particular,
\begin{align}
\frac{\mu(n)}{\log^2y}(1*\Lambda*\Lambda)(n)&=\mu(n)\Bigg(a_{20}'\Bigg(\frac{\log n}{\log y}\Bigg)^{2}+a_{21}'\sum_{p\mid n}\Bigg(\frac{\log p}{\log y}\Bigg)^2\Bigg),
\end{align}
\begin{align}
\frac{\mu(n)}{\log^2y}(1*\Lambda*\Lambda*\Lambda)(n)&=\mu(n)\Bigg(a_{30}' \Bigg(\frac{\log n}{\log y}\Bigg)^{3}+a_{32}'\Bigg(\frac{\log n}{\log y}\Bigg)\sum_{p\mid n}\Bigg(\frac{\log p}{\log y}\Bigg)^2\Bigg)+a_{33}'\sum_{p\mid n}\Bigg(\frac{\log p}{\log y}\Bigg)^3\Bigg).
\end{align}
Substituting these two formulae into \eqref{+1.18}, we have
\begin{align}
\label{x+1.28}
a(n)=\mu(n)\bigg(P_1\bigg(\frac{\log y/n}{\log y}\bigg)+\sum_{2\le k\le K}P_k\bigg(\frac{\log y/n}{\log y}\bigg)\sum_{p\mid n}\bigg(\frac{\log p}{\log y}\bigg)^k\bigg)
\end{align}
for $K\le3$. We take the following simple expression
\begin{align}
\label{x+1.29}
a(n)=\mu(n)\bigg(P_1\bigg(\frac{\log y/n}{\log y}\bigg)+P_2\bigg(\frac{\log y/n}{\log y}\bigg)\sum_{p\mid n}P\bigg(\frac{\log p}{\log y}\bigg)\bigg)
\end{align}
for convenience.

This short version of Feng's coefficient contains only one prime variable, which does not carry along some of difficult problems in Feng \cite{Feng}, and thus it makes the calculation simple. Also, we can easily see from this short version that Feng's coefficient works more effectively than Conrey's due to the present of the terms $\sum_{p\mid n}\big(\frac{\log p}{\log y}\big)^k$ with $k\ge2$. These terms can not be well approximated by the polynomial $P_1\big(\frac{\log y/n}{\log y}\big)$.

With the coefficient in \eqref{x+1.29}, we will face the same problem as in \cite{Feng} and \cite{PR} when estimate error terms for $\theta<\frac47$. But (B) of Theorem \ref{thm1} means that error terms can be controlled if we kick out large primes $p\ge y^{\frac34}$ in the second term. This dropping of large primes may bring a loss in final numeric results, but we may expect that it will be a weak loss. Let us see it in detail. For large $p$ with $p\ge y^{\frac34}$, the sum $\sum_{p\mid n, p\ge y^{\frac34}}$ has only one term, and  the difference between $\frac{\log p}{\log y}$ and $\frac{\log n}{\log y}$ is a multiple in $[\frac34, 1]$. This means that
\begin{align}
P_2\bigg(\frac{\log y/n}{\log y}\bigg)\sum_{p\mid n,~p\ge y^{\frac34}}P\bigg(\frac{\log p}{\log y}\bigg)
\end{align}
should be approximated `well' by $P_1\big(\frac{\log y/n}{\log y}\big)$ in $a(n)$. In view of this, we kick out large primes in the second term to have
\begin{align}\label{x1.32}
a(n)=\mu(n)\bigg(P_1\bigg(\frac{\log y/n}{\log y}\bigg)+P_2\bigg(\frac{\log y/n}{\log y}\bigg)\sum_{p\mid n,~p\le y^{3/4}}P\bigg(\frac{\log p}{\log y}\bigg)\bigg).
\end{align}
Also, our final numeric calculation verifies this expectation.\\

\noindent
{\bf Remark.} It should be noted that, similar results as those in Theorem \ref{thm2}, which are for the Riemann zeta-function only, are also obtained by Pratt, Robles, Zaharescu and Zeindler \cite{PRZZ}. They worked independently of us, and their work was presented in arXiv soon after ours. Unlike our simplification of coefficients by ignoring non-squarefree $n$ in the beginning, they reserve these $n$ to keep inside $\Lambda$ in their coefficient to the last. This turns out to be a very smart action in the estimation of error terms. Reserving non-squarefree terms causes the calculation of the main term to be more complex, but a return comes when one treats error terms. Without the squarefree condition on $n$ in error terms, a Vaughan's identity generated by Heath-Brown in \cite{HeaV} then can be used to split $n$ well, and so the error terms can be well controlled by Weil's bound and Deshouillers and Iwaniec's estimate on Kloosterman sums without kicking out large prime factors. They also prove in detail that these non-squarefree $n$ terms in the main term contribute an error, thus one may reserve non-squarefree $n$ terms in the error terms only and then remove the condition $p\le y^\frac34$ in \eqref{x1.32} with the help of \cite{PRZZ}.

\section{Notation and some standard results}
As usual, we use $\epsilon$ to denote an arbitrarily small positive constant that may vary from line to line in the following.

Both $\exp(x)$ and $e^x$ denote the exponential function, but $e(x)$ means $\exp(2\pi i x)$.

Throughout the paper  $q$ is a positive integer, which may go to infinity as $T\rightarrow \infty$, and $\chi$ is a character to the modulus $q$. Gauss' sum is
\begin{align}
\tau(\chi)=\sum_{n=1}^q\chi(n)e(n/q).
\end{align}
It is known that $\tau(\chi)\le q^{1/2}$, and in particular
\begin{align}
\tau(\chi)\tau(\overline\chi)=\chi(-1)q
\end{align}
when $\chi$ is a primitive character. Also we will use the following well-known formula
\begin{align}\label{++2.3}
\sum_{j=1}^q\chi(j)e\bigg(\frac{nj}{q}\bigg)=\overline{\chi}(n)\tau(\chi).
\end{align}

A Dirichlet $L$-function $L(s,\chi)$ has the following functional equation
\begin{align}
\label{FE}
 L(1-s,\chi)=H(1-s,\chi)L(s,\overline{\chi}),
\end{align}
where
\begin{align}\label{+2.4}
   H(1-s,\chi)=(2\pi)^{-s}q^{s-1}\tau(\chi)\Gamma(s)\bigg(e^{-\pi is/2}+\chi(-1)e^{\pi is/2}\bigg).
\end{align}
Due to this functional equation, $L(s,\chi)$ has an analytic continuation to the whole plane with a possible pole at $s=1$, which occurs with residue $q^{-1}\phi(q)$ when $\chi$ is a principal character $\chi_0$ .

Let $0<x\le1$, the Hurwitz zeta-function is defined by
\begin{align}
\zeta(s,x)=\sum_{n\ge1}(n+x)^{-s}
\end{align}
for Re$(s)>1$. The functional equation for the Hurwitz zeta-function is
\begin{align}
\label{m13}
\zeta(1-s,x)=\frac{\Gamma(s)}{(2\pi)^s}\Bigg\{e^{-\pi is/2}F(s,x)+e^{\pi is/2}F(s,-x)\Bigg\},
\end{align}
where
\begin{align}\label{+2.7}
F(s,x)=\sum_{n=1}^\infty e(nx)n^{-s}
\end{align}
for Re$(s)>1$.

\section{The proposition and the proof of Theorem \ref{thm1}}
In this section, we present a proposition and then prove Theorem \ref{thm1} from this proposition as in \cite{BCH} and \cite{Con-More}.

\begin{proposition}
\label{pro}
Let $\chi$ be a primitive Dirichlet character (mod $q$) with $\log q=o(\log T)$ and $\alpha=a/\mathcal{L}_\chi,~\beta=b/\mathcal{L}_\chi$ with $a,~b\in\mathbb{C}$ and $a,~b\ll1$. Suppose that $\eta>0$, $\Delta=T^{1-\eta}$, $y=T^\theta$ and $s_0=1/2+iw$ with $T\le w\le2T$. Let
\begin{align}\label{x3.1}
g(\alpha,\beta,w)=\frac1{i\Delta\pi^{1/2}}\int\limits_{(1/2)}e^{(s-s_0)^2 \Delta^{-2}}L(s+\alpha,\chi) L(1-s+\beta,\overline{\chi})B(s,\chi)B_1(1-s,\overline{\chi})ds
\end{align}
with $(c)$ denoting the straight line path from $c-i\infty$ to $c+\infty$ and where
\begin{align}
B_1(1-s,\overline{\chi})=\sum_{n\le y}\frac{\overline{\chi(n)}\overline{a(n)}}{n^{1-s}}.
\end{align}
Then we have, uniformly in $a,~b,$ and $w$, that
\begin{align}\label{x3.3}
g(\alpha,\beta,w)=\sum_{h,k\le y}&\frac{\chi_0(hk)(h,k)^{\alpha+\beta}a(h)\overline{a(k)}}{h^{1+\beta} k^{1+\alpha}}\notag\\
&\times \Bigg(L(1-\alpha-\beta,\chi_0)\bigg(\frac {2\pi hk}{qw(h,k)^2}\bigg)^{\alpha+\beta} +L(1+\alpha+\beta,\chi_0)\Bigg)+O(T^{-\epsilon_\theta}),
\end{align}
and in the particular case $\alpha=\beta=0$,
\begin{align}\label{x3.4}
g(w)=\frac{\phi(q)}{q}\sum_{h,k\le y}\frac{a(h)\overline{a(k)}}{[h,k]}\chi_0(hk)\bigg(\log\frac{wq(h,k)^2}{2\pi hk}+2\gamma+c_q\bigg)+O(T^{-\epsilon_\theta})
\end{align}
with $c_q=\sum_{p\mid q}(\log p)/(p-1)$. Here $\epsilon_\theta$ is a constant depending on $\theta$ in the following ways:
\begin{description}
  \item[(A1)] In general, $\epsilon_\theta>0$ for any given $\theta<\frac{17}{33}$;
  \item[(B1)] When $a(n)$ has a special form as in (B) of Theorem \ref{thm1}, $\epsilon_\theta>0$ for any given $\theta<\frac{4}{7}$.
\end{description}

\end{proposition}

We prove Theorem \ref{thm1} from this proposition exactly as in Section 3 of Balasubramanian, Conrey and Heath-Brown\cite{BCH} and Section 5 of Conrey\cite{Con-More}. Or more specifically, let
\begin{align}
g_Q(\alpha,\beta,w)=\frac1{i\Delta\pi^{1/2}}\int\limits_{(1/2)}e^{(s-s_0)^2\Delta^{-2}}V(s+\alpha,\chi) V(1-s+\beta,\overline{\chi})B(s,\chi)B_1(1-s,\overline{\chi})ds.
\end{align}
Using the proposition as in Section 5 of Conrey\cite{Con-More} we have
\begin{align}\label{x3.8}
&g_Q(\alpha,\beta,w)=Q\bigg(\frac{-d}{da}\bigg)\overline{Q}\bigg(\frac{-d}{db}\bigg) \\
&\Bigg(\sum_{h,k\le y}\frac{\chi_0(hk)(h,k)^{\alpha+\beta}a(h)\overline{a(k)}}{h^{1+\beta} k^{1+\alpha}}\Bigg(L(1-\alpha-\beta,\chi_0)\bigg(\frac {2\pi hk}{qw(h,k)^2}\bigg)^{\alpha+\beta} +L(1+\alpha+\beta,\chi_0)\Bigg)\Bigg)+O(T^{-\epsilon_\theta})\notag
\end{align}
with $\alpha=a/\mathcal{L}_\chi,~\beta=b/\mathcal{L}_\chi$. Then it follows exactly as in Section 3 of Balasubramanian, Conrey and Heath-Brown\cite{BCH} that
\begin{align}
I_R(Q,\chi)=\int_T^{2T}g_Q(\alpha,\beta,w)|_{a=b=-R}dw+O(T^{1-\epsilon_\theta}),
\end{align}
which gives \eqref{+1.10} with the help of \eqref{x3.8}. Similarly, we can obtain \eqref{+1.12} from \eqref{x3.4}, and thus we prove Theorem \ref{thm1}.\\

\noindent
{\bf Remark.} It should be reminded that formula \eqref{x3.8} is obtained directly by adding $Q\big(\frac{-d}{da}\big)$ and $\overline{Q}\big(\frac{-d}{db}\big)$ to the main term of $g(\alpha,\beta,w)$ without any action on the error term. Since $g$ is analytic in the complex variables $a$ and $b$ for $a,~b\ll1$, Cauchy's integral formula enables it. This action will also arise in Section \ref{sec5} when we calculate the main term, and it will not be pointed out any more.

\section{the main term of the proposition}
In this section, we produce the main term of the proposition after preparing some pivo-tal lemmas.
\begin{lemma}
\label{lemmx4.1}
Suppose that $1<c<2$, then
\begin{align}\label{x4.1}
\frac1{i\Delta\pi^{1/2}}\int\limits_{(c)}e^{(s-z)^2\Delta^{-2}}\Gamma(s)(2\pi i x)^{-s}ds=\int\limits_0^\infty v^{z}\exp\bigg(-\frac{\Delta^2\log^2v}{4}\bigg)e(-xv)\frac{dv}{v}
\end{align}
for any $x\neq0,~z$ and $\Delta>0$.
\end{lemma}
This lemma was exploited in the proof of Lemma 2 \cite{BCH}, which proved \eqref{x4.1} by the theory of Mellin transforms.

\begin{lemma}
\label{lemm1}
Suppose that $1<c<2$ and $q,~\Delta,~s_0,~\beta$ are given as in Proposition \ref{pro}. Let
\begin{align}
J(x,s_0,\beta,\Delta,\chi)=\frac1{i\Delta\pi^{1/2}}\int\limits_{(c)}e^{(s-s_0)^2\Delta^{-2}}H(1-s+\beta,\chi)x^{-s}ds
\end{align}
with $H(1-s+\beta,\chi)$ defined by \eqref{+2.4}.
Then
\begin{align}
J=\frac{\tau(\chi)}{qx^{\beta}}\int\limits_0^\infty v^{s_0-\beta}\exp\bigg(-\frac{\Delta^2\log^2v}{4}\bigg)\Big(e\Big(-\frac{xv}{q}\Big) +\chi(-1)e\Big(\frac{xv}{q}\Big)\Big)\frac{dv}{v}
\end{align}
for any $x\neq0$.
\end{lemma}
\begin{proof}
By a variable change $s-\beta\rightarrow s$, we have
\begin{align}
\label{m4}
J=\frac{\tau(\chi)}{qx^\beta}(J_1+\chi(-1)J_2),
\end{align}
where
\begin{align}
J_1=\frac1{i\Delta\pi^{1/2}}\int\limits_{(c)}e^{(s+\beta-s_0)^2\Delta^{-2}}\Gamma(s)\bigg(2\pi e^{\frac{\pi i}{2}}\frac xq\bigg)^{-s}ds
\end{align}
and
\begin{align}
J_2=\frac1{i\Delta\pi^{1/2}}\int\limits_{(c)}e^{(s+\beta-s_0)^2\Delta^{-2}}\Gamma(s)\bigg(2\pi e^{\frac{-\pi i}{2}}\frac xq\bigg)^{-s}ds.
\end{align}
Then we deduce $J_1$ and $J_2$ directly by Lemma \ref{lemmx4.1} and prove the lemma.
\end{proof}

\begin{lemma}
\label{lemm2}Let $H,~K$ be integers $(K\ge1)$ such that any two of $q,H,K$ are coprime. Suppose that $\alpha,~\beta,~s\in\mathbb{C}$ and let
\begin{align}
D\bigg(s,\alpha,\beta,\frac{H}{Kq},\chi\bigg)=\sum_{m,n}\frac{\chi(m)\chi(n)}{m^{s+\alpha}n^{s+\beta}}e\bigg(\frac{mnH}{Kq} \bigg),
\end{align}
then
\begin{align}
&D\bigg(s,\alpha,\beta,\frac{H}{Kq},\chi\bigg)-K^{1-2s-\alpha-\beta}\tau(\chi)\chi(K)\overline{\chi}(H)\notag\\
&\times\bigg(q^{-s-\alpha}L(s+\beta,\chi_0)\zeta(s+\alpha)+q^{-s-\beta}L(s+\alpha,\chi_0) \zeta(s+\beta)-q^{-2s-\alpha-\beta}\phi(q)\zeta(s+\alpha)\zeta(s+\beta)\bigg)\notag
\end{align}
is an entire function of $s$. Also, $D$ satisfies the equation
\begin{align}
D\bigg(1-s,\alpha,\beta,&\frac{H}{Kq},\chi\bigg)=\frac2{(Kq)^2}\bigg(\frac{Kq}{2\pi}\bigg)^{2s-\alpha-\beta} \Gamma(s-\alpha) \Gamma(s-\beta)\notag\\
&\times\bigg\{\cos\frac\pi2 (2s-\alpha-\beta)A_1\bigg(s,\alpha,\beta,\frac{H}{Kq},\chi\bigg)+\cos\frac\pi2 (\alpha-\beta)A_2\bigg(s,\alpha,\beta,\frac{H}{Kq},\chi\bigg)\bigg\}
\end{align}
with
\begin{align}
A_1\bigg(s,\alpha,\beta,\frac{H}{Kq},\chi\bigg)=\sum_{1\le v,u\le Kq}\chi(u)\chi(v)e\Big(\frac{uvH}{Kq}\Big)F(s-\alpha,\frac{u}{Kq})F(s-\beta,\frac{v}{Kq})
\end{align}
and
\begin{align}
A_2\bigg(s,\alpha,\beta,\frac{H}{Kq},\chi\bigg)=\sum_{1\le v,u\le Kq}\chi(u)\chi(v)e\Big(\frac{uvH}{Kq}\Big)F(s-\alpha,\frac{u}{Kq})F(s-\beta,-\frac{v}{Kq}).
\end{align}
Here $F(s,x)$ is defined in \eqref{+2.7}. Moreover, we have $D\Big(0,\alpha,\beta,\frac{H}{Kq},\chi\Big) \ll_\epsilon q^{3/2+\epsilon}K^{1+\epsilon}$ for any $\epsilon>0$ when $\alpha,\beta\ll(\log Kq)^{-1}$.
\end{lemma}

\begin{proof}
The entire function $D$ has been pointed out in \cite{Bau}, and one can deduce the upper bound for $D\Big(0,\alpha,\beta,\frac{H}{Kq},\chi\Big)$ directly from this entire function. Thus, we only prove the functional equation here. It follows from the definition of $D$ that
\begin{align}
\label{m12}
 D\bigg(s,\alpha,\beta,\frac{H}{Kq},\chi\bigg)=(Kq)^{-2s-\alpha-\beta}\sum_{1\le v,u\le Kq}\chi(u)\chi(v)e\Big(\frac{uvH}{Kq}\Big)\zeta\bigg(s+\alpha,\frac{u}{Kq}\bigg)\zeta\bigg(s+\beta,\frac{v}{Kq}\bigg)
\end{align}
for Re$(s)\ge1-\min\{\text{Re}(\alpha),\text{Re}(\beta)\}$. By analytic continuation, one may obtain that (\ref{m12}) is available in the whole plane. Making the variable change $s\rightarrow 1-s$ and employing the functional equation of the Hurwitz zeta-function (\ref{m13}) in (\ref{m12}), we have
\begin{align}\label{+3.14}
 D\bigg(1-s,\alpha,\beta,\frac{H}{Kq},\chi\bigg)=&\frac1{(Kq)^2}\bigg(\frac{Kq}{2\pi}\bigg)^{2s-\alpha-\beta} \Gamma(s-\alpha) \Gamma(s-\beta)\sum_{1\le v,u\le Kq}\chi(u)\chi(v)e\bigg(\frac{uvH}{Kq}\bigg)\notag\\
 &\times\Bigg\{e^{-\pi i(s-\alpha)/2}F\bigg(s-\alpha,\frac{u}{Kq}\bigg)+e^{\pi i(s-\alpha)/2}F\bigg(s-\alpha,-\frac{u}{Kq}\bigg)\Bigg\}\notag\\
 &\times\Bigg\{e^{-\pi i(s-\beta)/2}F\bigg(s-\beta,\frac{v}{Kq}\bigg)+e^{\pi i(s-\beta)/2}F\bigg(s-\beta,-\frac{v}{Kq}\bigg)\Bigg\}.
\end{align}
It is easy to see following two formulae
\begin{align}
&\sum_{1\le v,u\le Kq}\chi(u)\chi(v)e\bigg(\frac{uvH}{Kq}\bigg)\Bigg\{F\bigg(s-\alpha,\frac{ u}{Kq}\bigg)F\bigg(s-\beta,\frac{ v}{Kq}\bigg) -F\bigg(s-\alpha,\frac{-u}{Kq}\bigg)F\bigg(s-\beta,\frac{- v}{Kq}\bigg)\Bigg\}=0
\end{align}
and
\begin{align}
&\sum_{1\le v,u\le Kq}\chi(u)\chi(v)e\bigg(\frac{uvH}{Kq}\bigg)\Bigg\{F\bigg(s-\alpha,\frac{ u}{Kq}\bigg)F\bigg(s-\beta,\frac{-v}{Kq}\bigg) -F\bigg(s-\alpha,\frac{-u}{Kq}\bigg)F\bigg(s-\beta,\frac{v}{Kq}\bigg)\Bigg\}=0.
\end{align}
Then the lemma follows when we expand the two brackets in \eqref{+3.14} and use these two formulae to simplify it.
\end{proof}

\begin{lemma}
\label{lemm3}
Let $H,~K$ be integers with $K\ge1$ and any two of $q,~H,~K$ are coprime. Suppose that $\alpha,~\beta,~x\in\mathbb{C}$ with $\alpha\neq\beta$, $\alpha,~\beta\neq1$, Im$(x)>0$ and let
\begin{align}
S\bigg(x,\alpha,\beta,\frac{H}{Kq},\chi\bigg)=\sum_{m,n}\frac{\chi(m)\chi(n)}{m^{\alpha}n^{\beta}} e\bigg(\frac{mnH}{Kq}\bigg)e(mnx).
\end{align}
Then, for any $c>1-\min\{\text{Re}(\alpha),\text{Re}(\beta)\}$,
\begin{align}
S\bigg(x,\alpha,\beta,&\frac{H}{Kq},\chi\bigg)\notag\\
=&L(1-\alpha+\beta,\chi_0)K^{-1+\alpha-\beta}q^{-1}\tau(\chi)\overline{\chi}(H) \chi(K)z^{-1+\alpha}\Gamma(1-\alpha)\notag\\
&+L(1-\beta+\alpha,\chi_0)K^{-1+\beta-\alpha}q^{-1}\tau(\chi)\overline{\chi}(H) \chi(K)z^{-1+\beta}\Gamma(1-\beta)\notag\\
&+D(0,\alpha,\beta,\frac{H}{Kq},\chi)+\frac1{(Kq)^2\pi i}\int\limits_{(c)}z^{s-1}\Gamma(1-s)\Gamma(s-\alpha) \Gamma(s-\beta)\bigg(\frac{Kq}{2\pi}\bigg)^{2s-\alpha-\beta}\notag\\
&\times\bigg\{\cos\frac\pi2 (2s-\alpha-\beta)A_1\bigg(s,\alpha,\beta,\frac{H}{Kq},\chi\bigg)+\cos\frac\pi2 (\alpha-\beta)A_2\bigg(s,\alpha,\beta,\frac{H}{Kq},\chi\bigg)\bigg\}ds
\end{align}
with $z=-2\pi ix$.
\end{lemma}
\begin{proof}
By Mellin's formula,
\begin{align}
S&=\sum_{m,n}\frac{\chi(m)\chi(n)}{m^{\alpha}n^{\beta}} e\bigg(\frac{mnH}{Kq}\bigg)\frac1{2\pi i}\int\limits_{(c)}\Gamma(s)(-2\pi imnx)^{-s}ds\notag\\
&=\frac1{2\pi i}\int\limits_{(c)}D\bigg(s,\alpha,\beta,\frac{H}{Kq},\chi\bigg)\Gamma(s)z^{-s}ds,
\end{align}
where $c$ could be any real number that satisfies $c>1-\min\{\text{Re}(\alpha),\text{Re}(\beta)\}$. We move the path of integration to $(1-c)$ and cross three simple poles at $1-\alpha$, $1-\beta$ and $0$ for $\alpha\neq\beta$ and $\alpha,~\beta\neq1$. Since $L(s,\chi_0)$ has residue $q^{-1}\phi(q)$ at $s=1$, we have by Lemma \ref{lemm2} that residues of these poles are equal to
\begin{align}
&L(1-\alpha+\beta,\chi_0)K^{-1+\alpha-\beta}q^{-1}\tau(\chi)\overline{\chi}(H) \chi(K)z^{-1+\alpha}\Gamma(1-\alpha)\notag\\
&+L(1-\beta+\alpha,\chi_0)K^{-1+\beta-\alpha}q^{-1}\tau(\chi)\overline{\chi}(H) \chi(K)z^{-1+\beta}\Gamma(1-\beta)+D(0,\alpha,\beta,\frac{H}{Kq},\chi).
\end{align}
If we make the change of variable $s\rightarrow 1-s$ and use the functional equation in Lemma \ref{lemm2}, one will see that the integration on $(1-c)$ evolves into
\begin{align}
&\frac1{(Kq)^2\pi i}\int\limits_{(c)}z^{s-1}\Gamma(1-s)\Gamma(s-\alpha) \Gamma(s-\beta)\bigg(\frac{Kq}{2\pi}\bigg)^{2s-\alpha-\beta}\notag\\
&\times\bigg\{\cos\frac\pi2 (2s-\alpha-\beta)A_1\bigg(s,\alpha,\beta,\frac{H}{Kq},\chi\bigg)+\cos\frac\pi2 (\alpha-\beta)A_2\bigg(s,\alpha,\beta,\frac{H}{Kq},\chi\bigg)\bigg\}ds.
\end{align}
Thus the lemma follows by Cauchy's theorem.
\end{proof}

\begin{lemma}
\label{lemmx4.5}
Let $\alpha,~\beta\ll\log^{-1}T$, $0<\delta<\pi/2$, $z=\frac12+\beta+i\omega$ with $T\le\omega\le2T$ and $\Delta=T^{1-\eta}$ with $\eta>0$. We define
\begin{align}
\label{x4.20}
r_\delta(z,\alpha)=\int\limits_{L_\delta}v^{z}\exp\bigg(-\frac{\Delta^2\log^2v} {4}\bigg)(v-1)^{-1+\alpha}\frac{dv}{v},
\end{align}
where $L_{\delta}$ is the half-line $L_{\delta}=\{re^{i\delta}:r>0\}$. Let
\begin{align}
\label{x4.21}
W(z,\alpha)=\Gamma(1-\alpha)\{(-2\pi i)^\alpha r_{\delta}(z,\alpha)-(2\pi i)^\alpha r_{-\delta}(z,\alpha)\},
\end{align}
then we have
\begin{align}
W(z,\alpha)=-2\pi i\bigg(\frac{w}{2\pi}\bigg)^{-\alpha}+O(T^{-\eta}).
\end{align}
\end{lemma}

\begin{proof}
We denote $W$ by
\begin{align}
W(z,\alpha)=W_1(z,\alpha)-W_2(z,\alpha)
\end{align}
with obvious meanings. We consider $W_1$ first, and $W_2$ will be treated similarly. By definition,
\begin{align}
\label{x4.30}
W_1(z,\alpha)=-2\pi i\int\limits_{L_{\delta}}v^{z}\exp\bigg(-\frac{\Delta^2\log^2v}{4}\bigg) \Gamma(1-\alpha)(-2\pi i(v-1))^{-1+\alpha}\frac{dv}{v}.
\end{align}
Let $\theta(v)=\arg(-2\pi i(v-1))$, then it is easy to see that $-\frac\pi2+\delta<\theta(v)<\frac\pi2$ for $\delta<\arg(v-1)<\pi$ with $v$ in $L_\delta$. It is well-known that
\begin{align}
\Gamma(s)=\int_0^\infty u^{s-1}e^{-u}du,
\end{align}
and then we have
\begin{align}
\Gamma(1-\alpha)(-2\pi i(v-1))^{-1+\alpha}&=\int_0^\infty u^{-\alpha}e^{-u}(-2\pi i(v-1))^{-1+\alpha}du\notag\\
&=\int\limits_{L_{-\theta(v)}} u^{-\alpha}e(u(v-1))du
\end{align}
by making the change of variable $u\rightarrow-2\pi i(v-1)u$ in the integral. It is easy to check that
\begin{align}
\int\limits_{L_{\theta}} u^{-\alpha}e(u(v-1))du
\end{align}
is absolutely convergent for any $\theta$ with $|\theta+\theta(v)|<\frac\pi2$,
then Cauchy's theorem tells us that the integral can be moved to any path $L_\theta$ when $|\theta+\theta(v)|<\frac\pi2$. Thus we move the path to $L_{-\frac\delta2}$ and have
\begin{align}
\Gamma(1-\alpha)(-2\pi i(v-1))^{-1+\alpha}=\int\limits_{L_{-\frac\delta2}} u^{-\alpha}e(u(v-1))du.
\end{align}
When we apply this formula in \eqref{x4.30} and interchange the order of integrals, we have
\begin{align}
W_1(z,\alpha)=-2\pi i\int\limits_{L_{-\frac\delta2}} u^{-\alpha}e(-u)du\int\limits_{L_{\delta}}v^{z}\exp\bigg(-\frac{\Delta^2\log^2v} {4}\bigg) e(uv)\frac{dv}{v}.
\end{align}
We calculate the last integral by Lemma \ref{lemmx4.1}, and then
\begin{align}
W_1(z,\alpha)=-2\pi i
\frac1{i\Delta\pi^{1/2}}\int\limits_{(c)}e^{(s-z)^2\Delta^{-2}}\Gamma(s)ds \int\limits_{L_{-\frac\delta2}} u^{-\alpha}e(-u)(-2\pi i u)^{-s}du
\end{align}
with $1<c<2$.
It is easy to verify that the last integral is absolutely convergent.
We make the change of variable $u\rightarrow (-2\pi i)^{-1}u$ in the second integral to have
\begin{align}
W_1(z,\alpha)&=(-2\pi i)^\alpha
\frac1{i\Delta\pi^{1/2}}\int\limits_{(c)}e^{(s-z)^2\Delta^{-2}}\Gamma(s)ds \int\limits_{L_{-\frac\pi2-\frac\delta2}}e^uu^{-s-\alpha}du\notag\\
&=(-2\pi i)^\alpha
\frac1{i\Delta\pi^{1/2}}\int\limits_{(c)}e^{(s-z)^2\Delta^{-2}}\Gamma(s)ds \int\limits_{L_{-\pi}}e^uu^{-s-\alpha}du
\end{align}
by moving the last integral to the path $L_{-\pi}$. By the same way, one may obtain a similar expression for $W_2$
\begin{align}
W_2(z,\alpha)=(2\pi i)^\alpha
\frac1{i\Delta\pi^{1/2}}\int\limits_{(c)}e^{(s-z)^2\Delta^{-2}}\Gamma(s)ds \int\limits_{L_{\pi}}e^uu^{-s-\alpha}du.
\end{align}
Thus
\begin{align}
W(z,\alpha)=\frac{(2\pi)^\alpha}{i\Delta\pi^{1/2}}\int\limits_{(c)}e^{(s-z)^2 \Delta^{-2}}\Gamma(s) \Bigg(e^{-\frac\pi 2 i\alpha}\int\limits_{L_{-\pi}}e^uu^{-s-\alpha}du-e^{\frac\pi 2 i\alpha}\int\limits_{L_{\pi}}e^uu^{-s-\alpha}du\Bigg)ds.
\end{align}
We make variable changes $u\rightarrow e^{-i\pi}u$ and $u\rightarrow e^{i\pi}u$ in the two integrals on $L_{-\pi}$ and $L_{\pi}$ respectively, then
\begin{align}
e^{-\frac\pi 2 i\alpha}\int\limits_{L_{-\pi}}e^uu^{-s-\alpha}du-e^{\frac\pi 2 i\alpha}\int\limits_{L_{\pi}}e^uu^{-s-\alpha}du&=(e^{-\pi i(1-s-\frac\alpha2)}-e^{\pi i(1-s-\frac\alpha2)})\int_0^\infty e^{-u}u^{-s-\alpha}du\notag\\
&=-2i\sin(\pi(1-s-\frac\alpha2))\Gamma(1-s-\alpha).
\end{align}
It is well-known that
\begin{align}
\sin(\pi(1-s-\frac\alpha2))=\frac\pi{\Gamma(s+\frac\alpha2) \Gamma(1-s-\frac\alpha2)}.
\end{align}
Thus we have
\begin{align}
W(z,\alpha)&=-2\pi i\frac{(2\pi)^{\alpha}}{i\Delta\pi^{1/2}}\int\limits_{(c)}e^{(s-z)^2 \Delta^{-2}}\frac{\Gamma(s)\Gamma(1-s-\alpha)}{\Gamma(s+\frac\alpha2) \Gamma(1-s-\frac\alpha2)}ds\notag\\
&=-2\pi i\frac{(2\pi)^{\alpha}}{i\Delta\pi^{1/2}}\int\limits_{(\frac12)}e^{(s+\beta-z)^2 \Delta^{-2}}\frac{\Gamma(s+\beta)\Gamma(1-s-\beta-\alpha)}{\Gamma(s+\beta+\frac\alpha2) \Gamma(1-s-\beta-\frac\alpha2)}ds+O(\exp(-T^{2\eta}))
\end{align}
by making the variable change $s\rightarrow s+\beta$ and then moving the integral to $(\frac12)$. In moving the path we cross a simple pole at $1-\beta-\alpha$ with residue $\ll\exp(-T^{2\eta})$.
The Stirling's approximation gives
\begin{align}
\frac{\Gamma(s+\beta)\Gamma(1-s-\beta-\alpha)}{\Gamma(s+\beta+\frac\alpha2) \Gamma(1-s-\beta-\frac\alpha2)}=t^{-\alpha}(1+O(t^{-1})),
\end{align}
then it follows that
\begin{align}
\label{x4.38}
W(z,\alpha)=-2\pi i\frac1{\Delta\pi^{\frac12}}\int\limits_{-\infty}^{\infty}e^{-(t-w)^2 \Delta^{-2}}\bigg(\frac{t}{2\pi}\bigg)^{-\alpha} dt+O(T^{-1}).
\end{align}
Let $A$ be a large constant, we note that $(\frac{t}{2\pi})^{-\alpha}\ll 1$ for $t\ll T^A$ and $(\frac{t}{2\pi})^{-\alpha}\ll e^{\epsilon (t-w)^2\Delta^{-2}}$ for any $\epsilon>0$. Thus the integral in \eqref{x4.38} on $\{t:|t-w|\ge 2\Delta \log T\}$ is actually an error term $\ll\exp(-\log^2 T)$. For $|t-w|\le2 \Delta\log T$, it is easy to see that
\begin{align}
\bigg(\frac{t}{2\pi}\bigg)^{-\alpha}&= \bigg(\frac{w}{2\pi}\bigg)^{-\alpha} \bigg(1+O\bigg(\alpha\frac{t-w}{w}\bigg)\bigg) =\bigg(\frac{w}{2\pi}\bigg)^{-\alpha}(1+O(T^{-\eta})).
\end{align}
Thus we have
\begin{align}
W(z,\alpha)=-2\pi i\bigg(\frac{w}{2\pi}\bigg)^{-\alpha}+O(T^{-\eta})
\end{align}
and this completes our proof of the lemma.
\end{proof}

We now begin to prove the proposition. We split $g(\alpha,\beta,w)$ into the main term and error terms. In this section we produce the main term by following the method of \cite{Con-More}, and the error terms will be estimated in next section.

Let $\eta_0>0$ be a small and fixed real number.  We move the path of integration in \eqref{x3.1} to $(c)$ with $c=1+\eta_0$. It is easy to see that $|\alpha|,~|\beta|<\eta_0$ for sufficiently large $T$ since $\alpha,~\beta\ll1/\mathcal{L}_\chi$. Thus, if $\chi$ is a principal Dirichlet character, in moving the path we cross a pole at $s=1-\alpha$. However, the contribution of this residue is negligible since
\begin{align}
\exp((1-\alpha-s_0)^2\Delta^{-2})\ll\exp(-T^{2\eta}),
\end{align}
which decays rapidly as $T\rightarrow\infty$. By the functional equation (\ref{FE}) for $L(1-s+\beta,\chi)$ and Lemma \ref{lemm1}, we can interchange summation and integration to have
\begin{align}
g(\alpha,\beta,w)=&\sum_{h,k\le y}\frac{\chi(h)\overline{\chi}(k)a(h)\overline{a(k)}}{k} \sum_{m,n}\chi(m)\chi(n)m^{-\alpha}n^{\beta} J\bigg(\frac{mnh}{k},s_0,\beta,\Delta,\chi\bigg)+O(\exp(-T^{2\eta}))\notag\\
=&\frac{\tau(\overline{\chi})}{q}\sum_{h,k\le y}\frac{\chi(h)\overline{\chi}(k)a(h)\overline{a(k)}}{h^{\beta}k^{1-\beta}} \sum_{m,n}\chi(m)\chi(n)m^{-\alpha-\beta} \int\limits_0^\infty v^{s_0-\beta}\exp\bigg(-\frac{\Delta^2\log^2v}{4}\bigg)\notag\\
&\times\bigg(\overline{\chi}(-1)e\Big(\frac{mnhv}{kq}\Big) +e\Big(-\frac{mnhv}{kq}\Big)\bigg)\frac{dv}{v}+O(\exp(-T^{2\eta})).
\end{align}

Let $\delta>0$ be a small real number and $L_\delta$ be the half-line given in Lemma \ref{lemmx4.5}. We express the integral as a sum of two integrals and use Cauchy's theorem to move one path to $L_\delta$ and the other to $L_{-\delta}$. We interchange the integration and the summation over $m,~n$ to have
\begin{align}
\label{m28}
g(\alpha,\beta,w)=\frac{\tau(\overline{\chi})}{q}\sum_{h,k\le y}\frac{\chi(h)\overline{\chi}(k)a(h)\overline{a(k)}}{h^{\beta}k^{1-\beta}} (\overline{\chi}(-1)I_1+I_2)+O(\exp(-T^{2\eta})),
\end{align}
where
\begin{align}
I_1=\int\limits_{L_\delta}v^{s_1}\exp\bigg(-\frac{\Delta^2\log^2v}{4}\bigg) S\bigg(\frac{H(v-1)}{Kq},\alpha+\beta,0,\frac{H}{Kq},\chi\bigg)\frac{dv}{v}
\end{align}
and
\begin{align}
I_2=\int\limits_{L_{-\delta}}v^{s_1}\exp\bigg(-\frac{\Delta^2\log^2v}{4}\bigg) S\bigg(-\frac{H(v-1)}{Kq},\alpha+\beta,0,-\frac{H}{Kq},\chi\bigg)\frac{dv}{v}
\end{align}
with $s_1=s_0-\beta$ and $S$  defined as in Lemma \ref{lemm3}. Here $H=h/(h,k)$ and $K=k/(h,k)$. Then by Lemma \ref{lemm3}, we have
\begin{align}
I_1=M_1+R_1+E_1,
\end{align}
where
\begin{align}
M_1=&\frac{\tau(\chi)\chi(K)\overline{\chi}(H)}{q}\int\limits_{L_{\delta}}v^{s_1} \exp\bigg(-\frac{\Delta^2\log^2v}{4}\bigg) \notag\\ &\times\Bigg[L(1-\alpha-\beta,\chi_0)\Gamma(1-\alpha-\beta)K^{-1+\alpha+\beta}\bigg(-2\pi i\frac{H}{Kq}(v-1)\bigg)^{-1+\alpha+\beta}\notag\\
&+L(1+\alpha+\beta,\chi_0)K^{-1-\alpha-\beta}\bigg(-2\pi i\frac{H}{Kq}(v-1)\bigg)^{-1}\Bigg]\frac{dv}{v},
\end{align}
\begin{align}
R_1=D\bigg(0,\alpha+\beta,0,\frac{H}{Kq},\chi\bigg)\int\limits_{L_{\delta}}v^{s_1} \exp\bigg(-\frac{\Delta^2\log^2v}{4}\bigg)\frac{dv}{v}
\end{align}
and
\begin{align}
E_1=\int\limits_{L_{\delta}}v^{s_1} \exp\bigg(-\frac{\Delta^2\log^2v}{4}\bigg)F_1(v)\frac{dv}{v}
\end{align}
with
\begin{align}
F_1(v)&=\frac1{(Kq)^2\pi i}\int\limits_{(c)}\bigg(-2\pi i\frac{H}{Kq}(v-1)\bigg)^{s-1}\Gamma(1-s)\Gamma(s-\alpha-\beta) \Gamma(s)\bigg(\frac{Kq}{2\pi}\bigg)^{2s-\alpha-\beta}\notag\\
&\times\bigg\{\cos\frac\pi2 (2s-\alpha-\beta)A_1\bigg(s,\alpha+\beta,0,\frac{H}{Kq},\chi\bigg)+\cos\frac\pi2 (\alpha+\beta)A_2\bigg(s,\alpha+\beta,0,\frac{H}{Kq},\chi\bigg)\bigg\}ds.
\end{align}
There are similar expressions for $I_2=M_2+R_2+E_2$.

Now we come to deduce the main term of the proposition, which will come from $M_1$ and $M_2$. Let $r_\delta(s_1,\alpha)$ be defined as in Lemma \ref{lemmx4.5}, we have
\begin{align}
M_1=\frac{\tau(\chi)\chi(K)\overline{\chi}(H)}{q}\bigg[&L(1-\alpha-\beta,\chi_0) \Gamma(1-\alpha-\beta)\bigg(-2\pi i\frac Hq\bigg)^{-1+\alpha+\beta}r_\delta(s_1,\alpha+\beta)\notag\\
&+L(1+\alpha+\beta,\chi_0)\bigg(-2\pi i\frac Hq\bigg)^{-1}K^{-\alpha-\beta}r_\delta(s_1,0)\bigg]
\end{align}
and
\begin{align}
M_2=\frac{\tau(\chi)\chi(K)\overline{\chi}(-H)}{q}\bigg[&L(1-\alpha-\beta,\chi_0) \Gamma(1-\alpha-\beta)\bigg(2\pi i\frac Hq\bigg)^{-1+\alpha+\beta}r_{-\delta}(s_1,\alpha+\beta)\notag\\
&+L(1+\alpha+\beta,\chi_0)\bigg(2\pi i\frac Hq\bigg)^{-1}K^{-\alpha-\beta}r_{-\delta}(s_1,0)\bigg].
\end{align}
We use $M_1$ and $M_2$ to substitute $I_1$ and $I_2$ in (\ref{m28}) and obtain the main term of $g(\alpha,\beta,w)$
\begin{align}
\label{x4.53}
(-2\pi i)^{-1}\sum_{h,k\le y}&\frac{\chi_0(hk)(h,k)^{\alpha+\beta}a(h)\overline{a(k)}}{h^{1+\beta} k^{1+\alpha}}\notag\\
&\times\Bigg(L(1-\alpha-\beta,\chi_0)\bigg(\frac {HK}q\bigg)^{\alpha+\beta}W(s_1,\alpha+\beta) +L(1+\alpha+\beta,\chi_0)W(s_1,0)\Bigg)
\end{align}
with $W(s_1,\cdot)$ given by \eqref{x4.21}. Here we have used the fact that $\tau(\chi)\tau(\overline{\chi})=\chi(-1)q$ for a primitive character $\chi$. Then by Lemma \ref{lemmx4.5} we have
\begin{align}
W(s_1,0)=-2\pi i+O(T^{-\eta})
\end{align}
and
\begin{align}
\label{x4.55}
W(s_1,\alpha+\beta)=-2\pi i\bigg(\frac{w}{2\pi}\bigg)^{-\alpha-\beta}+O(T^{-\eta}).
\end{align}
We substitute these two formulae into \eqref{x4.53} and get the main term of $g(\alpha,\beta,w)$ in the proposition.

Since
\begin{align}
L(1+s,\chi_0)=\frac{\phi(q)}{q}\bigg(\frac1s+\gamma+c_q+o(1)\bigg)
\end{align}
as $s\rightarrow0$, then \eqref{x3.4} follows directly from \eqref{x3.3} with $\alpha,~\beta\rightarrow0$. Thus we get the main term of $g(w)$ in the proposition.

\section{The error terms}
In this section we consider the error terms arising in the above section. Since $\chi(-1)$ equal $1$ or $-1$, we actually need to bound
\begin{align}
\label{e1}
\frac{\tau(\overline{\chi})}{q}\sum_{h,k\le y}\frac{\chi(h)\overline{\chi}(k)a(h)a(k)}{h^{\beta}k^{1-\beta}}(R_i+E_i)
\end{align}
for $i=1$ and $2$. It is not difficult to note that these two situations are identical, so we deal with $i=1$ only. For $R_1$, by Lemma 4 of \cite{BCH}
\begin{align}
\label{e2}
\int\limits_{L_\delta}v^{s_1}\exp\bigg(-\frac{\Delta^2\log^2v}{4}\bigg)\frac{dv}{v}\ll\exp(-T^{2\eta}),
\end{align}
then we see from the upper bound $D\Big(0,\alpha,\beta,\frac{H}{Kq},\chi\Big)\ll_\epsilon q^{\frac32+\epsilon}K^{1+\epsilon}$ given in Lemma \ref{lemm2} that
\begin{align}
R_1\ll T^{-20}.
\end{align}
This means that the contribution to (\ref{e1}) from $R_1$ is negligible.

From the definition of $E_1$, we may split the part of (\ref{e1}) which involves $E_1$ into two terms, one of which is
\begin{align}
\label{e4}
Z=\frac{\tau(\overline{\chi})}{q}\int\limits_{L_\delta}\int\limits_{(c)}G(\alpha+\beta,v,s_1,\Delta,s) \mathscr{M}(\alpha,\beta,s,\chi)dsdv,
\end{align}
where
\begin{align}\label{defG}
G(\alpha,v,s_1,\Delta,s)=&\frac1{2\pi^2} v^{s_1}\exp\bigg(-\Delta^2\frac{\log^2v}{4}\bigg)\Gamma(s)\Gamma(1-s)\Gamma(s-\alpha) (2\pi)^{\alpha-s}\notag\\
&\times\cos(\pi/2(2s-\alpha))e^{-\pi is/2}(v-1)^{s-1}v^{-1}
\end{align}
and
\begin{align}
\label{e6}
\mathscr{M}(\alpha,\beta,s,\chi)=q^{s-1-\alpha-\beta}\sum_{h,k\le y}\frac{\chi(h)\overline{\chi}(k)a(h)a(k)}{H^{1-s+\beta}K^{2-s+\alpha}(h,k)} A_1\bigg(s,\alpha+\beta,0,\frac{H}{Kq},\chi\bigg).
\end{align}
Here $H=h/(h,k)$ and $K=k/(h,k)$ as before, and $A_1$ is defined in Lemma \ref{lemm2}. The other term of $E_1$ may be treated in the same way as this one will be.

To estimate $Z$ we should firstly take out $A_1$ in the right-hand side of \eqref{e6}. Recalling the definition of $A_1$ we have that
\begin{align}
\label{e8}
A_1\bigg(s,\alpha+\beta,0,\frac{H}{Kq},\chi\bigg)=&\sum_{1\le v,u\le Kq}\chi(u)\chi(v)e\bigg(\frac{uvH}{Kq}\bigg)F\bigg(s-\alpha-\beta,\frac{u}{Kq}\bigg)F\bigg(s,\frac{v}{Kq}\bigg)\notag\\
=&\sum_{1\le v,u\le Kq}\chi(u)\chi(v)e\bigg(\frac{uvH}{Kq}\bigg)\sum_{m}\frac1{m^{s-\alpha-\beta}}e\bigg(\frac{mu}{Kq}\bigg) \sum_{n}\frac1{n^{s}}e\bigg(\frac{nv}{Kq}\bigg).
\end{align}
 Let the least positive residue mod $Kq$ of $Hu$ be $a$ in (\ref{e8}) and suppose that $\overline{H}$ is defined by
\begin{align}
H\overline{H}\equiv1~(\text{mod}~Kq), ~\ \ \ \ 0<\overline{H}\le Kq.
\end{align}
Then we have
\begin{align}
\label{e9}
A_1\bigg(s,\alpha+\beta,0,\frac{H}{Kq},\chi\bigg)=&\sum_{1\le a,v\le Kq}\chi(a\overline{H})\chi(v)e\bigg(\frac{av}{Kq}\bigg)\sum_{m}\frac1{m^{s-\alpha-\beta}} e\bigg(\frac{ma\overline{H}}{Kq}\bigg)\sum_{n}\frac1{n^{s}}e\bigg(\frac{nv}{Kq}\bigg)\notag\\
=&\sum_{1\le a\le Kq}\chi(a\overline{H})\sum_{m}\frac1{m^{s-\alpha-\beta}}e\bigg(\frac{ma\overline{H}}{Kq}\bigg) \sum_{n}\frac1{n^{s}}\sum_{1\le v\le Kq}\chi(v)e\bigg(\frac{(a+n)v}{Kq}\bigg).
\end{align}
We write $v=rq+j$ and have the sum over $v$ above equal to
\begin{align}\label{+5.10}
\sum_{{1\le j\le q}\atop{0\le r\le K-1}}\chi(rq+j)e\bigg(\frac{(a+n)(rq+j)}{Kq}\bigg)
=\sum_{1\le j\le q}\chi(j))e\bigg(\frac{(a+n)j}{Kq}\bigg)\sum_{0\le r\le K-1}e\bigg(\frac{(a+n)r}{K}\bigg).
\end{align}
We note that
\begin{align}\label{+5.11}
\sum_{0\le r\le K-1}e\bigg(\frac{(a+n)r}{K}\bigg)=\begin{cases}
K &n\equiv-a~(\text{mod} K),\\
0 &\text{otherwise}.
\end{cases}
\end{align}
Thus we restrict $(n+a)/K$ to be an integer.
Then by \eqref{++2.3}
we have
\begin{align}\label{+5.12}
\sum_{1\le j\le q}\chi(j)e\Big(\frac{(a+n)j}{Kq}\Big)=\overline{\chi}\Big(\frac{a+n}{K}\Big)\tau(\chi).
\end{align}
Employing \eqref{+5.10}-\eqref{+5.12} into (\ref{e9}) we have
\begin{align}
\label{e15}
A_1\bigg(s,\alpha+\beta,0,\frac{H}{Kq},\chi\bigg)=&\chi(\overline{H})\tau(\chi)K\sum_{1\le a\le Kq}\chi(a)\sum_{m}\frac1{m^{s-\alpha-\beta}}e\bigg(\frac{ma\overline{H}}{Kq}\bigg)\sum_{n\equiv -a (\text{mod}~K)}\frac1{n^{s}}\overline{\chi}\bigg(\frac{a+n}{K}\bigg)\notag\\
=&\chi(\overline{H})\tau(\chi)K\sum_{m}\frac1{m^{s-\alpha-\beta}}\sum_{n}\frac1{n^{s}}\sum_{{1\le a\le Kq}\atop{a\equiv -n (\text{mod}~K)}}\chi(a)e\bigg(\frac{ma\overline{H}}{Kq}\bigg)\overline{\chi}\bigg(\frac{a+n}{K}\bigg).
\end{align}
To the sum over $a$ above, we denote $a=-n+Kj$. One notes that $j$ exactly runs through all residue classes mod $q$ when $a$ takes all possible values in the sum. Also we note that every term in the sum over $a$ is independent from the exact value of $j$ but depends on its residue class mod $q$. Thus, we have
\begin{align}
\sum_{{1\le a\le Kq}\atop{a\equiv -n (\text{mod}~k)}}\chi(a)e\bigg(\frac{ma\overline{H}}{Kq}\bigg)\overline{\chi}\bigg(\frac{a+n}{K}\bigg)
=\sum_{j=0}^{q-1} \chi(-n+Kj)e\bigg(\frac{m(-n+Kj)\overline{H}}{Kq}\bigg)\overline{\chi}(j).
\end{align}
Employing this into (\ref{e15}) we get
\begin{align}
\label{e17}
A_1\bigg(s,\alpha+\beta,0,\frac{H}{Kq},\chi\bigg)=&\chi(\overline{H})\tau(\chi)K\sum_{j=0}^{q-1} \sum_{m}\frac1{m^{s-\alpha-\beta}}\sum_{n}\frac1{n^{s}} e\bigg(\frac{m(-n+Kj)\overline{H}}{Kq}\bigg)\overline{\chi}(j)\chi(-n+Kj).
\end{align}

We now come to $\mathscr{M}$. Employing (\ref{e17}) into (\ref{e6}) and arranging the sums over $h$ and $k$ according to $g\equiv(h,k)$ we get
\begin{align}
\label{c11}
\mathscr{M}(\alpha,\beta,s,\chi)=&q^{s-1-\alpha-\beta}\tau(\chi)\sum_{g\le y}\frac{\chi_0(g)}{g}\sum_{j=0}^{q-1}\sum_{m,n}\frac{1}{m^{s-\alpha-\beta}n^s}\notag\\
&\times\sum_{{H,K\le y/g}\atop{(H,K)=1}}\frac{\chi_0(H)a(Hg)a(Kg)}{H^{1-s+\beta}K^{1-s+\alpha}} e\bigg(\frac{m(-n+Kj)\overline{H}}{Kq}\bigg)\overline{\chi}(Kj)\chi(-n+Kj).
\end{align}
In the above formula we have used the fact that
\begin{align}
\chi(H)\chi(\overline{H})=\chi_0(H).
\end{align}
The right-hand side of (\ref{c11}) seems to be complicated to estimate. This will be alleviated if we can remove some influence of $\chi$.  When fixing $j$ and residue classes modulo $q$ of $m,~n$ and $K$, we  note that $\overline{\chi}(Kj)\chi(-n+Kj)$ is a constant. It should be noted that this fixing also splits the right-hand side of (\ref{c11}) into $q^4$ sums. Let $j_1,j_2,j_3,~j_4$ be four constants with $0\le j_i\le q-1$. We designate $j=j_1$ and residue classes modulo $p$ of $m,n,K$ to be $j_2,j_3,j_4$ respectively, then
\begin{align}
e\bigg(\frac{m(-n+Kj)\overline{H}}{Kq}\bigg)\overline{\chi}(Kj)\chi(-n+Kj)=e\bigg(\frac{-mn\overline{H}+Kj_1j_2 \overline{H}}{Kq}\bigg)\overline{\chi}(j_1j_4)\chi(j_1j_4-j_3),
\end{align}
while, on the right-hand side, the last two factors are independent from exact values of $m,~n,~H,~K$ and $\ll1$ now. It should be highlighted that the fixing does not contain the variable $H$, and this is important for our following discussion. Thus $\mathscr{M}$ is split into $q^4$ sums of the shape
\begin{align}
C\tau(\chi)\sum_{g\le y}\frac{\chi_0(g)}{g}\sum_{M,N,U,V}\mathscr{M}_1(M,N,U,V,\alpha,\beta,g,s),
\end{align}
where $C\ll1$ and
\begin{align}
\label{c15}
\mathscr{M}_1(M,N,U,V,\alpha,\beta,g,s)=&q^{s-1-\alpha-\beta}\mathop{{\sum}'}_{m\sim M}\mathop{{\sum}'}_{n\sim N}\frac{1}{m^{s-\alpha-\beta}n^s}\notag\\
&\times \sum_{{u\sim U}\atop{(u,q)=1}}\mathop{{\sum}'}_{{v\sim V}\atop{(u,v)=1}}\frac{a(ug)a(vg)}{u^{1-s+\beta}v^{1-s+\alpha}}e\bigg(\frac{-mn\overline{u}+vj_1j_2\overline{u}}{vq}\bigg),
\end{align}
where $\mathop{{\sum}'}$ denotes the sum of $m,n$ and $K$ over fixed residue classes $j_2,j_3,j_4$ respectively. Here, the notation $x\sim X$ means $X<x\le2X$, and the sums on $U$ and $V$ have $\log y$ terms with $U,V\ll y/g$ and the sums on $M,N$ are for $M=2^I, ~N=2^J$ with $I,J=0,1,2,\cdots$. When devoted  to obtaining an upper bound uniform in $g$ for $\mathscr{M}_1$, one may note that the sum on $g$ just contributes a multiple $\log y$. For convenience, we get rid of it by multiplying $T^\epsilon$ in the following. Employing these into (\ref{e4}) we have that $Z$ is a sum over $M, N, U,V$ of terms of the shape
\begin{align}\label{defZ1}
Z_1(M,N,U,V)\ll T^\epsilon\frac{|\tau(\chi)|^2}{q}\int\limits_{L_\delta}\int\limits_{(c)} G(\alpha+\beta,v,s_1,\Delta,s) \mathscr{M}_1(M,N,U,V,\alpha,\beta,g,s)dsdv.
\end{align}
Thus, it is enough to estimate these $Z_1$ to bound $Z$. We carry this out by classifying $Z_1$ into two cases:
\begin{description}
  \item[Case one] $MNT^{1-3\eta}>UV$;
  \item[Case two] $MNT^{1-3\eta}\le UV$.
\end{description}
Here $\eta>0$ is the small constant in the proposition.
We will see that $Z_1$ in case one can be controlled well even for very large $U$ and $V$, and the estimation in case two limits the upper bound of $\theta$. When we estimate $Z_1$, we will need the following lemma.
\begin{lemma}
\label{leme1}
Let $G$ be defined by \eqref{defG} with $s_1,~\Delta$, and $\alpha$ as above. Suppose that $\delta=1/T$ and $c_1\le c_2$ are any given constants. Then we have uniformly in $\alpha$, $s_1$ and $c$ with $c_1\le c\le c_2$ that
\begin{align}
\int\limits_{L_\delta}\int\limits_{(c)}(1+|s|)|G(\alpha+\beta,v,s_1,\Delta,s)|dsdv\ll_\epsilon\Delta^{-2c-\frac32}T^{c+\frac32+\epsilon} +\Delta^{-c-\frac52}T^{\frac52+\epsilon}\notag
\end{align}
for any $\epsilon>0$.
\end{lemma}

\begin{proof}
This lemma mainly extends the scope of $c$ to $c_1\le c\le c_2$ for any given $c_1$ and $c_2$, while only the case when $c$ is a little greater than $0$ or $1$ was the concern of Lemma 5 in \cite{BCH}. The proof of this lemma is almost the same as Lemma 5 in \cite{BCH}, thus we only point out the difference between them. Since $c_1\le c\le c_2$, the upper bound
\begin{align}
|\Gamma(s)|\ll(1+|t|)^{c-\frac12}e^{-\frac{\pi}2|t|}
\end{align}
by Stirling's formula is available. Let $v=xe^{i\delta}$. We have
\begin{align}
|(v-1)^{s-1}|=a_c(x,\delta)\exp(-tb(x,\delta)),
\end{align}
where
\begin{align}
a_c(x,\delta)=\big((x-1)^2+2x(1-\cos\delta)\big)^\frac{c-1}{2}
\end{align}
and
\begin{align}
b(x,\delta)=\arctan\frac{x\sin\delta}{x\cos\delta-1}.
\end{align}
The main difference is a new estimate of $a_c(x,\delta)$
that
\begin{equation}
a_c(x,\delta)\ll\left\{
\begin{aligned}
&\max\bigg\{1,\bigg(\frac{\log\Delta}{\Delta}\bigg)^{c-1}\bigg\} \ \ \ \ \ \ \ \ \ \ \ &\text{if} \ \ \ x\le1-\frac{\log\Delta}{\Delta};\\
&\bigg(\frac{\log\Delta}{\Delta}\bigg)^{c-1}+\delta^{c-1} \ \ \ \ &\text{if}\ \  |x-1|\le\frac{\log\Delta}{\Delta};\\
&\max\bigg\{x^{c-1},\bigg(\frac{\log\Delta}{\Delta}\bigg)^{c-1}\bigg\} \ \ \ \ &\text{if}\ \ \ x\ge1+\frac{\log\Delta}{\Delta}.
\end{aligned}
\right.
\end{equation}
This estimate can be obtained directly if we split $c$ in two situations $c-1\ge0$ and $c-1<0$. It is easy to see from \cite{BCH} that integrals over $x\le1-\frac{\log\Delta}{\Delta}$ and $x\ge1+\frac{\log\Delta}{\Delta}$ are due to $W_1,~W_2$ and $W_4$, and one can easily follow the method in \cite{BCH} to control them. For the integral on $1-\frac{\log\Delta}{\Delta}\le x\le1+\frac{\log\Delta}{\Delta}$, we have
\begin{align}
W_3\ll \frac{\log\Delta}{\Delta}\bigg(\frac{\log\Delta}{\Delta\delta}\bigg)^{c+\frac32} \bigg(\bigg(\frac{\log\Delta}{\Delta}\bigg)^{c-1}+\delta^{c-1}\bigg)\ll_\epsilon\Delta^{-2c-\frac32}T^{c+\frac32+\epsilon} +\Delta^{-c-\frac52}T^{\frac52+\epsilon}
\end{align}
for any $\epsilon>0$, which proves the lemma.
\end{proof}

\subsection{Error with large $M,N$}\label{sec4.1}
In this subsection, we estimate the contribution of $Z_1$ in case one $MNT^{1-3\eta}>UV$. For these $Z_1$, we move the $s$ path of the integration to $s=c+it$ with some constant $c>2$ to be specified later. In moving the path we cross poles at $s=w$ with $w=2,3,\cdots, [c]$. Since
\begin{align}
\mathscr{M}_1(M,N,U,V,\alpha,\beta,g,w)\ll q^w (MN)^{1-w+\epsilon}(UV)^{w+\epsilon},
\end{align}
the residue at $s=w$ is
\begin{align}
\frac1{2\pi^2}\Gamma(w)\Gamma(w-\alpha-\beta)&(2\pi)^{\alpha+\beta-w}\cos(\pi/2 (2w-\alpha-\beta))e^{-\pi iw/2}\mathscr{M}_1(M,N,U,V,\alpha,\beta,g,w)\notag\\
\times&\int\limits_{L_\delta}v^{s_1}(v-1)^{w-1}\exp\bigg(-\frac{\Delta^2\log^2v}{4}\bigg)\frac{dv}{v}\ll_w(MN)^{1-w+\epsilon}T^{-10}
\end{align}
for any $w\ge2$ when $U,V\ll T^A$ for some given constant $A>0$. The above integral is estimated as follows, with the binomial expansion of $(v-1)^w$ we express the integral into several terms, and each term $\ll\exp(-T^{2\eta})$ by \eqref{e2}. Since $M, N, U$ and $V$ take values in geometric progression, we sum over $M, N, U$ and $V$ to find that all these residues contribute an error $\ll_c T^{-10}\log^2T$ to $Z$.

On the new path it is trivial that
\begin{align}\label{4.29}
\mathscr{M}_1(M,N,U,V,\alpha,\beta,g,s)=\begin{cases}
O_\epsilon \bigg(q^{c-1}\Big(\frac{UV}{MN}\Big)^c (MN)^{1+\epsilon}\bigg) &\text{for any}~ \epsilon>0,\\
O_A \bigg(q^{c-1}\Big(\frac{UV}{MN}\Big)^c MN\bigg) &\text{with}~ MN\le T^{A}~ \text{for some constant A}.
\end{cases}
\end{align}
Employing this into \eqref{defZ1}, we have by Lemma \ref{leme1} that
\begin{align}
Z_1(M,N,U,V)\ll_\epsilon q^{c-1}\bigg(\frac{UV}{MN}\bigg)^c (MN)^{1+\epsilon}(\Delta^{-2c-\frac32}T^{c+\frac32+\epsilon} +\Delta^{-c-\frac52}T^{\frac52+\epsilon}).
\end{align}
We estimate the right-hand side of this formula in following two situations.

If $MN\ge (UV)^{\frac{c}{c-2}}$, ie $UV\le (MN)^\frac{c-2}{c}$, we have
\begin{align}
Z_1(M,N,U,V)&\ll_\epsilon (MN)^{-1+\epsilon}q^{c-1}\Big(\Delta^{-c}T^{(c+\frac32)\eta+\epsilon}+\Delta^{-c}T^{\frac52\eta+\epsilon}\Big)\notag\\
&\ll_\epsilon (MN)^{-1+\epsilon}\Delta^{-1}
\end{align}
for any $c\ge2$. Thus, summing $Z_1$ over $M,N,U,V$ we find that these $Z_1$ contribute an error $\ll\Delta^{-1+\epsilon}$ to $Z$.

If $MN\le (UV)^{\frac{c}{c-2}}$, then both sums on $M$ and $N$ have $\ll\log T$ terms, that is to say, $Z$ contains $\ll q^4\log^4 T$ terms of these $Z_1$. Remembering that $MNT^{1-3\eta}>UV$, we have from \eqref{4.29} that
\begin{align}
\mathscr{M}_1(M,N,U,V,\alpha,\beta,g,s)\ll q^{c-1}T^{c(1-3\eta)} MN.
\end{align}
Substituting this into \eqref{defZ1}, we have by Lemma \ref{leme1} that
\begin{align}
Z_1(M,N,U,V)&\ll_\epsilon q^{c-1}T^{c(1-3\eta)} MN(\Delta^{-2c-\frac32}T^{c+\frac32+\epsilon} +\Delta^{-c-\frac52}T^{\frac52+\epsilon})\notag\\
&\ll_\epsilon MN(T^{-c\eta+\frac32\eta+2\epsilon}+T^{-2c\eta+\frac52\eta+2\epsilon})
\end{align}
since $\log q=o(\log T)$. Taking $\epsilon=\frac14\eta$ we have
\begin{align}\label{4.34}
Z_1(M,N,U,V)&\ll_\eta MNT^{-(c-2)\eta}\ll_\eta (UV)^{\frac{c}{c-2}}T^{-(c-2)\eta}.
\end{align}
Let $U,~V\ll T^A$ for some given constant $A>0$. There is certainly a constant $c_{\eta,A}>2$ that the right-hand side of \eqref{4.34} $\ll T^{-1}$ for any $c\ge c_{\eta,A}$. With $c=c_{\eta,A}$ we multiply the quantity to find that these $Z_1$ contribute  an error $\ll T^{-1}\log^4T$ to $Z$.

Thus we conclude that if  $U,~V\ll T^A$ for some constant $A>0$, all $Z_1$ in case one $MNT^{1-3\eta}>UV$ contribute an error $\ll_\epsilon\Delta^{-1+\epsilon}$. That is to say, to prove the proposition it remains to estimate $Z_1$ in case two: $MNT^{1-3\eta}\le UV$. We will estimate these $Z_1$ for (A1) and (B1) of the proposition in two subsections respectively.

\subsection{Error with small $M,N$ for (A1)}
We estimate $Z_1$ in case two  $MNT^{1-3\eta}\le UV$ for (A1) in this subsection, then the first part of the proposition follows from this estimate with the help of the estimate in section \ref{sec4.1}.
It is obvious that both sums on $M$ and $N$ have $\ll\log T$ terms in this case. Thus the quantity of these $Z_1$ is $\ll q^4\log^4 T$. To estimate $Z_1$ we move the $s$ path of the integration to $s=1/2+it$, crossing a pole at $s=1$ with residue
\begin{align}
\frac1{2\pi^2}\Gamma(1-\alpha-\beta)&(2\pi)^{\alpha+\beta}\cos(\pi/2(2-\alpha-\beta))\mathscr{M}_1(M,N,U,V,\alpha,\beta,g,1)\notag\\
&\times\int\limits_{L_\delta}v^{s_1}\exp\bigg(-\frac{\Delta^2\log^2v}{4}\bigg)\frac{dv}{v}\ll_\eta T^{-10}.
\end{align}

Lemma \ref{leme1} indicates that an estimate of $Z_1$ in the new line would be obtained from an uniform upper bound of $\mathscr{M}_1$. For (A1), we deduce this uniform upper bound from the following lemma.
\begin{lemma} [Bettin and Chandee \cite{BC}] \label{lembc}
Let $\alpha_m,~\beta_u, \gamma_v$ be complex numbers, where $M\le m<2M$, $U\le u<2U$, and $V\le v<2V$. Then for any $\epsilon>0$, we have
\begin{align}
\sum_{m}\sum_{u}\sum_{(v,u)=1}\alpha_m\beta_u\gamma_ve\bigg(\frac{m\overline{u}}{v}\bigg)\ll&_\epsilon||\alpha||~||\beta||~||\gamma|| \bigg(1+\frac{M}{UV}\bigg)^\frac12\notag\\
&\times\bigg((MUV)^{\frac7{20}+\epsilon}(U+V)^{\frac14}+(MUV)^{\frac38+\epsilon}(MU+MV)^{\frac18}\bigg),\notag
\end{align}
where $||\cdot||$ denotes the $L_2$ norm.
\end{lemma}

Before using Lemma \ref{lembc}, we firstly adjust the expression of $\mathscr{M}_1$ in \eqref{c15} as follows
\begin{align}\label{4.37}
\mathscr{M}_1(M,N,U,V,\alpha,\beta,g,s)=&q^{s-1-\alpha-\beta}\mathop{{\sum}'}_{m\sim M}\mathop{{\sum}'}_{n\sim N}\frac{1}{m^{s-\alpha-\beta}n^s}\notag\\
&\times \sum_{{u\sim U}\atop{(u,q)=1}}\frac{a(ug)}{u^{1-s+\beta}}e\bigg(\frac{j_1j_2\overline{u}}{q}\bigg)\mathop{{\sum}'}_{{v\sim V}\atop{(u,v)=1}}\frac{a(vg)}{v^{1-s+\alpha}}e\bigg(\frac{-mn\overline{u}}{vq}\bigg).
\end{align}
Then we employ Lemma \ref{lembc} with variable changes $mn\rightarrow m,~u\rightarrow u$ and $vq\rightarrow v$ in the above equation. As $s=1/2+it$ and $MNT^{1-3\eta}\le UV$, we have
\begin{align}
\mathscr{M}_1(M,N,U,V,\alpha,\beta,g,s)&\ll_\epsilon(MNUV)^{\frac7{20}+\epsilon}(U+V)^{\frac14}+(MNUV)^{\frac38+\epsilon}(MNU+MNV)^{\frac18}\notag\\
&\ll_\epsilon y^{\frac{33}{20}}T^{-\frac7{20}+\frac{21}{20}\eta+\epsilon}+y^{\frac{15}{8}}T^{-\frac12+\frac32\eta+\epsilon}.
\end{align}
Combining this and Lemma \ref{leme1} with $c=1/2$ we have
\begin{align}
Z_1(M,N,U,V)&\ll_\epsilon y^{\frac{33}{20}}T^{-\frac{17}{20}+\frac{81}{20}\eta+\epsilon}+y^{\frac{15}{8}}T^{-1+\frac92\eta+\epsilon}.
\end{align}
By counting the number we have that these $Z_1$ contribute an error
\begin{align}
&\ll_\epsilon q^4\log^4T\Big( y^{\frac{33}{20}}T^{-\frac{17}{20}+\frac{81}{20}\eta+\epsilon}+y^{\frac{15}{8}}T^{-1+\frac92\eta+\epsilon}\Big)\notag\\
&\ll_\epsilon y^{\frac{33}{20}}T^{-\frac{17}{20}+\frac{81}{20}\eta+\epsilon}+y^{\frac{15}{8}}T^{-1+\frac92\eta+\epsilon}
\end{align}
to $Z$, which proves (A1) of the proposition with the help of the estimate in section \ref{sec4.1}.

\subsection{Error with small $M,N$ for (B1)}
We estimate $Z_1$ in case two $MNT^{1-3\eta}\le UV$ for (B1) and prove the second part of the proposition in this subsection. We move the $s$ path of the integration in $Z_1$ to $s=c+it$ with $c=\eta_0$, where $\eta_0>0$ is a small constant to be specified later. The residue at the pole $s=1$ can be neglected as before.
Thus we need to estimate $Z_1$ in the new path, which will be deduced from a sharp enough upper bound of $\mathscr{M}_1$. We will deduce this upper bound first, which necessitates much space to discuss, and then we obtain an estimate of $Z_1$ directly from this upper bound (see also the formula \eqref{4.99}) and Lemma \ref{leme1}.

Let us recall that
\begin{align}
\mathscr{M}_1(M,N,U,V,\alpha,\beta,g,s)=&q^{s-1-\alpha-\beta}\mathop{{\sum}'}_{m\sim M}\mathop{{\sum}'}_{n\sim N}\frac{1}{m^{s-\alpha-\beta}n^s}\notag\\
&\times \sum_{{u\sim U}\atop{(u,q)=1}}\mathop{{\sum}'}_{{v\sim V}\atop{(u,v)=1}}\frac{a(ug)a(vg)}{u^{1-s+\beta}v^{1-s+\alpha}}e\bigg(\frac{-mn\overline{u}+vj_1j_2\overline{u}}{vq}\bigg),
\end{align}
where
\begin{align}\label{4.43}
a(ug)=\mu(ug)\mathcal{F}(ug)+\mu(ug)\mathcal{F}_1(ug)(\mathcal{F}_2*\mathcal{F}_3)(ug).
\end{align}
According to \eqref{4.43}, we split $\mathscr{M}_1$ into two terms
\begin{align}\label{+4.44}
\mathscr{M}_1=\mathscr{M}_{10}+\mathscr{M}_{11}.
\end{align}
with obvious meanings. We will estimate $\mathscr{M}_{10}$ and $\mathscr{M}_{11}$ in next two subsections respectively.

\subsubsection{Estimate of $\mathscr{M}_{10}$}
Due to the separability of $\mathcal{F}_0$ and familiar properties of M\"{o}bius function we have
\begin{align}\label{+5.42}
\mathscr{M}_{10}\ll(MN)^{-c}(UVq)^{c-1}\sum_{m\sim M}\sum_{n\sim N}\sum_{v\sim V}r(m)r(n)r(v)\sum_{{u\sim U}\atop{(u,vqg)=1}}\mu(u)r^*(u)e\bigg(\frac{-mn\overline{u}+vj_1j_2\overline{u}}{vq}\bigg).
\end{align}
Here functions $r$ are the same as that in \cite{Con-More}, which may be different at each occurrence but all meet the following condition: $r(n)$ depends on its argument $n$ as well as $g, s, \alpha, \beta, j_i, q, M, N, U,$ and $V$, but it has $r(n)\ll_{\epsilon}n^{\epsilon}$ for any $\epsilon>0$ uniformly in all other arguments. In addition,
\begin{align}
r^*(u)=\mathcal{F}_0(u)u^{s-1-\beta}U^{1-c}
\end{align}
is also an $r$ function but smooth in its dependency on $u$, satisfying
\begin{align}
\frac{d}{du}r^*(u)\ll(1+|s|)u^{-1}r(u)
\end{align}
for some $r(u)$ and having the property of separability as $\mathcal{F}_0$. It is easy to see that $\mathcal{F}_0$ is an $r^*$ function and so are other $\mathcal{F}_i$ which meet $\mathcal{F}'_i(x)\ll\frac1x$. Also one may note that a product of two $r^*$ functions is also an $r^*$ function. In virtue of the functions $r$, the restriction that $m, n$ to be fixed residue classes modulo $q$ has been removed in \eqref{+5.42}. When denoting $mn$ by $n$ and $MN$ by $N$, we may combine the sums on $m$ and $n$ in $\mathscr{M}_{10}$ to have
\begin{align}
\label{c111}
\mathscr{M}_{10}\ll(N)^{-c}(UVq)^{c-1}|S|
\end{align}
with
\begin{align}
\label{c112}
S=\sum_{n\sim N}\sum_{v\sim V}r(n)r(v)\sum_{{u\sim U}\atop{(u,vqg)=1}}\mu(u)r^*(u)e\bigg(\frac{-n\overline{u}+vj_1j_2\overline{u}}{vq}\bigg).
\end{align}

Let $\Omega=U^{\frac14}$. It is easy to see that $u\ge U>\Omega$ for large $U$. Thus, when using the Vaughan's identity
\begin{align}
\zeta(s)^{-1}=\zeta(s)^{-1}\bigg(1-\zeta(s)\sum_{n\le\Omega}\mu(n)n^{-s}\bigg)^2+2\sum_{n\le \Omega}\mu(n)n^{-s}+\zeta(s)\bigg(\sum_{n\le\Omega}\mu(n)n^{-s}\bigg)^2
\end{align}
to split $\mu(u)$, we note that the second term in the right-hand of the above formula does not contribute anything. That is to say,
\begin{align}\label{4.47}
\mu(u)=c_1(u)+c_2(u),
\end{align}
where
\begin{align}\label{4.49}
c_1(u)=\sum_{{u_1u_2u_3=u}\atop{u_1>\Omega,~ u_2>\Omega}}\mu(u_3)c_{11}(u_1)c_{11}(u_2)
\end{align}
for some function $c_{11}(x)\ll x^\epsilon$, and
\begin{align}\label{4.50}
c_2(u)=-\sum_{{u_1u_2u_3=u}\atop{u_1\le\Omega,~u_2\le\Omega}}\mu(u_1)\mu(u_2).
\end{align}

This leads to $S=S_1+S_2$ with
\begin{align}
S_i=\sum_{n\sim N}\sum_{v\sim V}r(n)r(v)\sum_{{u\sim U}\atop{(u,vqg)=1}}c_i(u)r^*(u)e\bigg(\frac{-n\overline{u}+vj_1j_2\overline{u}}{vq}\bigg)
\end{align}
for $i=1,2$. One may note that $S_i$ here are similar to ones in \cite{Con-More} but with an extra factor $e\bigg(\frac{j_1j_2\overline{u}}{q}\bigg)$. We will treat this extra factor in different ways when estimating $S_1$ and $S_2$.

We come to $S_1$ first. We group together $u_3$ and the larger of $u_1$ and $u_2$ in \eqref{4.49} into a variable $b$ and name the other variable $a$. Then due to the separability of $r^*$ we split $S_1$ into $\ll_{\epsilon}y^\epsilon$ sums of the shape
\begin{align}
S_1'=\sum_{n\sim N}\sum_{v\sim V}r(n)r(v)\sum_{{{a\sim A}\atop{b\sim B}}\atop{(ab,v)=1}}r(a)r(b)e\bigg(\frac{j_1j_2\overline{a}\overline{b}}{q}\bigg)e\bigg(\frac{-n\overline{a}\overline{b}}{vq}\bigg),
\end{align}
where $U\ll AB\ll U$ and $\Omega\le A\le B$. We need to separate variables $a$ and $b$ in the coefficient, however the factor $e\Big(\frac{j_1j_2\overline{a}\overline{b}}{q}\Big)$ seems impossible to be separated. Actually, we do not separate it but dispel it in the following way. We note that the value of the factor $e\Big(\frac{j_1j_2\overline{a}\overline{b}}{q}\Big)$ depends on residue classes modulo $q$ of $a$ and $b$ but not their specific values. Thus, when we fix the residue classes modulo $q$ of $a$ and $b$ respectively, this factor is a constant. By classifying $a$ and $b$ according to residue classes modulo $q$ we split $S_1'$ into $q^2$ sums of shape
\begin{align}
S_1''=C\sum_{n\sim N}\sum_{v\sim V}r(n)r(v)\sum_{{{a\sim A}\atop{b\sim B}}\atop{(ab,v)=1}}r(a)r(b)e\bigg(\frac{-n\overline{a}\overline{b}}{vq}\bigg)
\end{align}
with constant $C\le1$ decided by residue classes of $a,~b$ modulo $q$. One notes that the restriction of $a, b$ to be fixed residue classes modulo $q$ in the sum has been removed since it can be absorbed by $r(a)$ and $r(b)$. It is easy to see that $S_1''$ is almost the same to $S_1'$ in \cite{Con-More} if we regard $vq$ here as $v$ in \cite{Con-More}. Since $\log q=o(\log T)$ and $UV\ge NT^{1-3\eta}$, we follow the way in \cite{Con-More} to obtain that
\begin{align}
S_1''\ll_\epsilon\max(TN,~UVq)(yN)^\epsilon T^{-\frac12+\epsilon}y^{\frac78}\ll_\epsilon UV T^{-\frac12+3\eta+2\epsilon}y^{\frac78}
\end{align}
by the following lemma.
\begin{lemma}\label{lem4.3}
Suppose that $V,~B,~N,~A\ge1$ and $|c(a,n)|\le1$. Then for any integer $d$ and $\epsilon>0$
\begin{align}
&\sum_{v\sim V}\sum_{{b\sim B}\atop{(b,v)=1}}\bigg|\sum_{n\sim N}\sum_{{a\sim A}\atop{(a,v)=1}}c(a,n)e\bigg(\frac{n\overline{d}\overline{a}\overline{b}}{v}\bigg)\bigg|\notag\\
&\ll(VBNA)^{\frac12+\epsilon}\{(VB)^{\frac12}+(N+A)^{\frac14}[VB(N+dA)(V+dA^2)+dNA^2B^2]^{\frac14}\}.\notag
\end{align}
\end{lemma}
This lemma is Lemma 1 of Deshouillers and Iwaniec \cite{DIP}. When treating $S_1''$, we actually use its special case with $d=1$, see also Lemma 9 in \cite{Con-More}.

Now we come to consider $S_2$. By \eqref{4.50} and the separability of $r^*$, we group together $u_1$ and $u_2$ into a variable $a$ and replace $u_3$ by $b$ to split $S_2$ into $\ll_\epsilon y^\epsilon$ sums of the shape
\begin{align}\label{+5.55}
S_2'=\sum_{v\sim V}r(v)\sum_{{b\sim B}\atop{(b,vg)=1}}r^*(b)\sum_{n\sim N}r(n)\sum_{{a\sim A}\atop{(a,v)=1}}r(a)e\bigg(\frac{-n\overline{a}\overline{b}+vj_1j_2\overline{a}\overline{b}}{vq}\bigg),
\end{align}
where $U\ll AB\ll U$ and $A\le \Omega^2=U^{\frac12}$. If $A\ge U^{\frac14}$, one notes that $S_2'$ reduces to the case $S_1$. If $A\le U^{\frac14}$, we can estimate the sum on $b$ by Weil's bound for the Kloosterman sum
\begin{align}\label{4.60}
\sum_{{b\sim B}\atop{(b,vg)=1}}e\bigg(\frac{l\overline{b}}{v}\bigg)\ll_{\epsilon}v^{1/2}(vg)^\epsilon (l,v)(1+Bv^{-1}).
\end{align}
Since
\begin{align}
\frac{d}{db}r^*(b)\ll(1+|s|)b^{-1}r(b),
\end{align}
we have by Abel's summation formula that
\begin{align}\label{+5.58}
\sum_{{b\sim B}\atop{(b,vqg)=1}}r^*(b) e\bigg(\frac{-n\overline{a}\overline{b}+vj_1j_2\overline{a}\overline{b}}{vq}\bigg)\ll(1+|s|) (yq)^\epsilon(Vq)^{1/2}((vj_1j_2-n)\overline{a},vq)\bigg(1+\frac{B}{vq}\bigg).
\end{align}
For $a\overline{a}\equiv1$ (mod $vq$) it follows that
\begin{align}\label{4.63}
((vj_1j_2-n)\overline{a},vq)&= ((vj_1j_2-n),vq)\le q((vj_1j_2-n),v)=q(n,v).
\end{align}
Thus we have by \eqref{+5.55}, \eqref{+5.58} and \eqref{4.63}
\begin{align}
S_2&\ll_\epsilon(1+|s|) (yNT)^\epsilon V^{\frac12}\bigg(1+\frac{B}{vq}\bigg)\sum_{n\sim N}\sum_{a\sim A}\sum_{v\sim V}(n,v)\notag\\
&\ll_\epsilon(1+|s|) (yNT)^\epsilon ANV^{\frac12}(V+B)\notag\\
&\ll_\epsilon(1+|s|) (yNT)^\epsilon (ANV^{\frac32}+NUV^{\frac12}).
\end{align}
It is obvious that
\begin{align}
ANV^{\frac32}+NUV^{\frac12}\ll Ny^{\frac74},
\end{align}
and this indicates
\begin{align}
S_2\ll_\epsilon(1+|s|)(TN)(yNT)^\epsilon T^{-1}y^{\frac74}\ll_\epsilon(1+|s|)UV T^{-1+3\eta+\epsilon}y^{\frac74}
\end{align}
as $UV\ge NT^{1-3\eta}$.
In conclusion, we have in any case
\begin{align}
S_2\ll_\epsilon(1+|s|) UV \Big(T^{-\frac12+3\eta+2\epsilon}y^{\frac78}+T^{-1+3\eta+\epsilon}y^{\frac74}\Big).
\end{align}

Thus we conclude from these upper bounds of $S_i$ for $i=1,2$ that
\begin{align}
S=S_1+S_2\ll_\epsilon(1+|s|) UVT^{\epsilon+3\eta} \Big(T^{-\frac12}y^{\frac78}+T^{-1}y^{\frac74}\Big)
\end{align}
for any $\epsilon>0$. Then it follows immediately by employing \eqref{c111} with $c=\eta_0$ that
\begin{align}\label{4.69}
\mathscr{M}_{10}&\ll(N)^{-\eta_0}(UVq)^{\eta_0-1}|S|\ll_\epsilon(1+|s|) \bigg(\frac{UV}{N}\bigg)^{\eta_0}T^{\epsilon+3\eta} \Big(T^{-\frac12}y^{\frac78}+T^{-1}y^{\frac74}\Big)\notag\\
&\ll_\epsilon(1+|s|) T^{\epsilon+3\eta+2\eta_0} \Big(T^{-\frac12}y^{\frac78}+T^{-1}y^{\frac74}\Big).
\end{align}

\subsubsection{Estimate of $\mathscr{M}_{11}$}
We now come to $\mathscr{M}_{11}$. It is easy to see that
\begin{align}
\mu(ug)(\mathcal{F}_2*\mathcal{F}_3)(ug)=0
\end{align}
for $(u,g)>1$. When $(u,g)=1$, due to the separability of $\mathcal{F}_2$ and $\mathcal{F}_3$, $(\mathcal{F}_2*\mathcal{F}_3)(ug)$ can be separated to finite terms of the form $(\mathcal{F}_2*\mathcal{F}_3)(g)\cdot(\mathcal{F}_2*\mathcal{F}_3)(u)$.
Since
\begin{align}
\mu(g)(\mathcal{F}_2*\mathcal{F}_3)(g)\ll g^\epsilon\ll T^\epsilon,
\end{align}
we have
\begin{align}\label{4.74}
\mathscr{M}_{11}\ll T^\epsilon(N)^{-c}(UVq)^{c-1}|\mathfrak{S}|,
\end{align}
where
\begin{align}
\label{4.75}
\mathfrak{S}=\sum_{n\sim N}\sum_{v\sim V}r(n)r(v)\sum_{{u\sim U}\atop{(u,vqg)=1}}r^*(u)\mu(u)(\mathcal{F}_2*\mathcal{F}_3)(u)e\bigg(\frac{-n\overline{u}+vj_1j_2\overline{u}}{vq}\bigg).
\end{align}
Here
\begin{align}
r^*(u)=\mathcal{F}_1(u)u^{s-1-\beta}U^{1-c},
\end{align}
which is also an $r^*$ function, and we also mark it by $r^*$. As in $\mathscr{M}_{10}$ we have combined sums on $m$ and $n$ to one variable, also denoted by $n$, in \eqref{4.75}.

When $U\le y^{\frac34}$, a trivial estimate gives that
\begin{align}
\mathfrak{S}\ll NUV\ll UVT^{-1+3\eta}y^{\frac74}
\end{align}
since $UV\ge NT^{1-3\eta}$ and $V\le y$.
Thus we always assume $U>y^{\frac34}$ in the following.

The difference between $\mathfrak{S}$ here and $S$ in $\mathscr{M}_{10}$ is an additional factor $(\mathcal{F}_2*\mathcal{F}_3)(u)$. Thus we need to separate $\mu(u)(\mathcal{F}_2*\mathcal{F}_3)(u)$ as $\mu(u)$ in $\mathscr{M}_{10}$. We note that
\begin{align}\label{+4.75}
\mu(u)(\mathcal{F}_2*\mathcal{F}_3)(u)&=\sum_{{u_1u_2=u}\atop{(u_1,u_2)=1}}\mu(u_1)\mu(u_2)\mathcal{F}_2(u_1) \mathcal{F}_3(u_2)\notag\\
&=\sum_{u_1u_2=u}\mu(u_1)\mu(u_2) \mathcal{F}_2(u_1)\mathcal{F}_3(u_2)\sum_{d\mid(u_1,u_2)}\mu(d)\notag\\
&=\sum_{d^2\mid u}\mu(d)\sum_{{u_1u_2=u/d^2}\atop{(u_1,d)=1,~(u_2,d)=1}}\mu(u_1)\mu(u_2)\mathcal{F}_2(du_1)\mathcal{F}_3(du_2).
\end{align}
Employing this with the properties of $\mathcal{F}_i$ into \eqref{4.75}, we have
\begin{align}
\mathfrak{S}\ll& T^\epsilon\sum_{d\sim D}\Bigg|\sum_{n\sim N}\sum_{v\sim V}r(n)r(v)\notag\\
&\times\sum_{{u_1\sim U_1}\atop{(u_1,vqdg)=1}}\sum_{{u_2\sim U_2}\atop{(u_2,vqdg)=1}}\mu(u_1)\mathcal{F}_2(u_1)r^*(u_1)\mu(u_2)r^*(u_2)e\bigg(\frac{-n\overline{d^2u_1u_2} +vj_1j_2\overline{d^2u_1u_2}}{vq}\bigg)\Bigg|
\end{align}
with $U\ll D^2U_1U_2\ll U$. Then we estimate $\mathfrak{S}$ by classifying them in the following three cases:
\begin{itemize}
  \item $DU_1\gg Uy^{-\frac34}$ and $DU_2\gg Uy^{-\frac34}$;
  \item $DU_1\ll Uy^{-\frac34}$;
  \item $DU_1\ge y^{\frac34}$.
\end{itemize}
It is easy to see these three cases contain all possible values of $D,~U_1$ and $U_2$ with $U\ll D^2U_1U_2\ll U$.

We start with the first case, which we denote by $\mathfrak{S}'$. We can also remove the factor $e\bigg(\frac{j_1j_2\overline{d^2u_1u_2}}{q}\bigg)$ by classifying $d,u_1$ and $u_2$ according to their residue classes modulo $q$, this splits $\mathfrak{S}'$ into no more than $q^3$ terms. It follows
\begin{align}\label{+5.74}
\mathfrak{S}'\ll T^\epsilon\sum_{d\sim D}\sum_{v\sim V}\sum_{{u_1\sim U_1}\atop{(u_1,vd)=1}}\bigg|\sum_{n\sim N}\sum_{{u_2\sim U_2}\atop{(u_2,vd)=1}}r(n)r(u_2)e\bigg(\frac{-n\overline{d^2u_1u_2}}{vq}\bigg)\bigg|.
\end{align}
This formula also holds if one interchanges $u_1$ and $u_2$ in \eqref{+5.74}. Thus we only need to treat the case $U_2\le U_1$ here. By Lemma \ref{lem4.3} we have
\begin{align}
\mathfrak{S}'&\ll T^\epsilon\sum_{d\sim D}(NVU_1U_2)^{\frac12}\Big\{(VU_1)^{\frac12}+(N+U_2)^{\frac14}[VU_1(N+d^2U_2)(V+d^2U_2^2)+d^2NU_1^2U_2^2]^{\frac14}\Big\}\notag\\
&\ll T^\epsilon\sum_{d\sim D}\bigg(\sum_{(\alpha_1,\alpha_2,\alpha_3,\alpha_4,\alpha_5)\in E} d^{\alpha_1}N^{\alpha_2} V^{\alpha_3} U_1^{\alpha_4} U_2^{\alpha_5}\bigg){^\frac14},
\end{align}
where
\begin{align}
E=\{(0,&2,4,4,2),(0,4,4,3,2),(2,4,3,3,4),(2,3,4,3,3),\notag\\
&(4,3,3,3,5),(2,2,4,3,4),(4,2,3,3,6)\}.
\end{align}
Let $A=dU_2$, we have $Uy^{-\frac34}\ll A\ll U^{\frac12}$ since $d^2U_1U_2\sim U$. It is easy to see
\begin{align}
\alpha_4+\alpha_5-\alpha_1\ge4
\end{align}
for all possible values that happened in $E$. Thus we have
\begin{align}
\mathfrak{S}'&\ll T^\epsilon\sum_{d\sim D}d^{-1}\bigg(\sum_{(a,n,u,v)\in E'} A^{a}N^{n} U^{u} V^{v}\bigg){^\frac14}\ll T^\epsilon\bigg(\sum_{(a,n,u,v)\in E'} A^{a}N^{n} U^{u} V^{v}\bigg){^\frac14}
\end{align}
with
\begin{align}
E'=\{(-2,&2,4,4),(-1,4,3,4),(1,4,3,3),(1,2,3,4),\notag\\
&(3,2,3,3),(2,3,3,3),(0,3,3,4)\}.
\end{align}
The set $E'$ is the same as $E$ given by formula (90) in \cite{Con-More}. One may note that the last two elements of them are different, it is possible due to a misprint in \cite{Con-More}. By employing $N\le UVT^{-1+3\eta}$, $A\ll U^{\frac12}$ and $U,V\le y$, we can treat every term with $a\ge0$ directly and obtain that
\begin{align}
A^{a}N^{n} U^{u} V^{v}\ll U^4V^4T^{-2+12\eta}y^{\frac72}
\end{align}
for these terms in $E'$ since $y\le T^{\frac47}$. For $a<0$, employing $N\le UVT^{-1+3\eta}$,~ $A\gg Uy^{-\frac34}$ and $U,V\le y$ we have
\begin{align}
A^{-2}N^{2} U^{4} V^{4}\ll U^4V^4T^{-2+6\eta}y^{\frac72}
\end{align}
and
\begin{align}
A^{-1}N^{4} U^{3} V^{4}\ll U^4V^4T^{-4+12\eta}y^{\frac{27}{4}}.
\end{align}
Thus we conclude that
\begin{align}\label{4.87}
\mathfrak{S}'\ll_\epsilon UVT^{-\frac12+\epsilon+3\eta}y^{\frac78}.
\end{align}

We now consider $\mathfrak{S}$ with $DU_1\ll Uy^{-\frac34}$, i.e. $DU_2\gg y^{\frac34}$. Let $\Omega'=Uy^{-\frac34}D^{-1}$, it is obvious that $1\le U_1\ll\Omega'\le y^{\frac14}D^{-1}<U_2\le u_2$. Thus, we split the M\"{o}bius function into two terms as in \eqref{4.47} and have
\begin{align}
\mu(u_2)=c_1'(u_2)+c_2'(u_2),
\end{align}
where $c_i'(n)$ are defined as $c_i(n)$ but with $\Omega$ replaced by $\Omega'$. This leads to $\mathfrak{S}=\mathfrak{S}_1+\mathfrak{S}_2$, where
\begin{align}
\mathfrak{S}_i\ll T^\epsilon\sum_{d\sim D}\bigg|\sum_{n\sim N}\sum_{v\sim V}r(n)r(v)\sum_{{u_1\sim U_1}\atop{(u_1,vqdg)=1}}\sum_{{u_2\sim U_2}\atop{(u_2,vqdg)=1}}r(u_1)c'_i(u_2)r^*(u_2)e\bigg(\frac{-n\overline{d^2u_1u_2} +vj_1j_2\overline{d^2u_1u_2}}{vq}\bigg)\bigg|.
\end{align}

To treat $\mathfrak{S}_1$, let us recall that
\begin{align}
c'_1(u_2)=\sum_{{u_{21}u_{22}u_{23}=u_2}\atop{u_{21}>\Omega', u_{22}>\Omega'}}\mu(u_{23})c'_{11}(u_{21})c'_{11}(u_{22})
\end{align}
with  some function $c_{11}'(x)\ll\frac1{x^\epsilon}$. Grouping together $u_1$, $u_{23}$ and the larger of $u_{21}$ and $u_{22}$ into a variable $b$ and renaming the other variable as $a$, we see that $\mathfrak{S}_1$ is split into $\ll_\epsilon y^\epsilon$ sums of the shape
\begin{align}
\mathfrak{S}_1'\ll T^\epsilon\sum_{d\sim D}\bigg|\sum_{n\sim N}\sum_{v\sim V}r(n)r(v)\sum_{{a\sim A}\atop{(a,vqdg)=1}}\sum_{{b\sim B}\atop{(b,vqdg)=1}}r(a)r(b)e\bigg(\frac{-n\overline{d^2u_1u_2}+vj_1j_2\overline{d^2u_1u_2}}{vq}\bigg)\bigg|,
\end{align}
where $U/D^2\ll AB\ll U/D^2$ and $Uy^{-\frac34}D^{-1}\ll A\le B$. Due to the fact that both $A$ and $B$ $\gg Uy^{-\frac34}D^{-1}$, we can treat $\mathfrak{S}_1'$ the same as $\mathfrak{S}'$ above. With \eqref{4.87} this means that
\begin{align}
\mathfrak{S}_1'\ll_\epsilon UVT^{-\frac12+\epsilon+3\eta}y^{\frac78}.
\end{align}

The treatment of $\mathfrak{S}_2$ is based on Weil's bound. Let us remember that
\begin{align}
c'_2(u_2)=-\sum_{{u_{21}u_{22}u_{23}=u_2}\atop{u_{21}\le\Omega', u_{22}\le\Omega'}}\mu(u_{21})\mu(u_{22}).
\end{align}
When grouping together $u_{21}$, $u_{22}$ and $u_1$ into a variable $a$ and replacing $u_{23}$ by $b$, we see that $\mathfrak{S}_2$ is split into $\ll_\epsilon y^\epsilon$ sums of the shape
\begin{align}
\mathfrak{S}_2'\ll T^\epsilon\sum_{d\sim D}\bigg|\sum_{n\sim N}r(n)\sum_{v\sim V}r(v)\sum_{{a\sim A}\atop{(a,vq)=1}}r(a)\sum_{{b\sim B}\atop{(b,vqdg)=1}}r^*(b)e\bigg(\frac{-n\overline{d^2ab}+vj_1j_2\overline{d^2ab}}{vq}\bigg)\bigg|,
\end{align}
where $U/D^2\ll AB\ll U/D^2$ and $A\le \Omega'^3$. It is obvious that that $\Omega'^4=U^4y^{-3}D^{-4}\le UD^{-2}\ll AB$, which means $B\gg \Omega'$. Thus, for $A\gg\Omega'$, $\mathfrak{S}_2'$ actually reduces to $\mathfrak{S}'$ above. If $A\ll\Omega'$, by using Weil's bound for the Kloosterman sum over $b$ and Abel's summation formula as before, we obtain that
\begin{align}
\mathfrak{S}_2'\ll_\epsilon (1+|s|)T^\epsilon V^{\frac12}\bigg(1+\frac{B}{Vq}\bigg)\sum_{d\sim D}\sum_{n\sim N}\sum_{a\sim A}\sum_{v\sim V}((vj_1j_2-n)\overline{d^2a},vq).
\end{align}
Similarly as \eqref{4.63} we have
\begin{align}
((vj_1j_2-n)\overline{d^2a},vq)\le q(n,v).
\end{align}
Since $\log q=o(\log T)$, it follows that
\begin{align}
\mathfrak{S}'_2&\ll_\epsilon (1+|s|)T^\epsilon V^{\frac12}\sum_{d\sim D}\bigg(1+\frac{B}{Vq}\bigg)\sum_{n\sim N}\sum_{a\sim A}\sum_{v\sim V}(n,v)\notag\\
&\ll_\epsilon(1+|s|) T^\epsilon \sum_{d\sim D}(ANV^{\frac12}(V+B)).
\end{align}
Due to $A\ll \Omega'\ll Uy^{-\frac34}d^{-1}$ and $U/d^2\ll AB\ll U/d^2$, we get
\begin{align}
\mathfrak{S}_2'&\ll_\epsilon(1+|s|) T^\epsilon \Big(y^{-\frac34}NUV^{\frac32}\sum_{d\sim D}d^{-1}+NUV^{\frac12}\sum_{d\sim D}d^{-2}\Big)\notag\\
&\ll_\epsilon(1+|s|)UV T^{-1+3\eta+\epsilon}y^{\frac74}
\end{align}
since $N\le UVT^{-1+3\eta}$ and $U,V\le y$.

It remains to treat $\mathfrak{S}$ with $DU_1\ge y^{\frac34}$. If $\mathcal{F}_2$ is separable in $\mathbb{F}$, we can treat it the same as the case $DU_1\ll Uy^{-\frac34}$. If  $\mathcal{F}_2$ is separable in $\Big\{\mathcal{F}:\mathcal{F}(x)\ll_\epsilon x^\epsilon,~\mathcal{F}(x)=0\ \ \text{for}\ \ x\ge y^{\frac34}\Big\}$, it vanishes since $\mathcal{F}_2(du_1)=0$ in \eqref{+4.75}.

We conclude that
\begin{align}
 \mathfrak{S}\ll_\epsilon(1+|s|)UVT^{-1+3\eta+\epsilon}y^{\frac74}+UVT^{-\frac12+3\eta+\epsilon}y^{\frac78}
\end{align}
in all cases.
Then it follows immediately by using \eqref{4.74} with $c=\eta_0$ that
\begin{align}
\mathscr{M}_{11}&\ll_\epsilon(N)^{-c}(UVq)^{c-1}|\mathfrak{S}|\ll_\epsilon(1+|s|) \bigg(\frac{UV}{N}\bigg)^{\eta_0}T^{\epsilon+3\eta} \Big(T^{-\frac12}y^{\frac78}+T^{-1}y^{\frac74}\Big)\notag\\
&\ll_\epsilon(1+|s|) T^{\epsilon+3\eta+2\eta_0} \Big(T^{-\frac12}y^{\frac78}+T^{-1}y^{\frac74}\Big).
\end{align}
With the help of  the estimate on $\mathscr{M}_{10}$ in last subsection, it follows that
\begin{align}\label{4.99}
\mathscr{M}_{1}\ll_\epsilon(1+|s|) T^{\epsilon+3\eta+2\eta_0} \Big(T^{-\frac12}y^{\frac78}+T^{-1}y^{\frac74}\Big).
\end{align}
Combining this and Lemma \ref{leme1} with $c=\eta_0$ we have
\begin{align}
Z_1&\ll_\epsilon T^{2\epsilon+3\eta+2\eta_0} \Big(T^{-\frac12}y^{\frac78}+T^{-1}y^{\frac74}\Big)\Big(\Delta^{-2\eta_0-\frac32}T^{\eta_0+\frac32+\epsilon} +\Delta^{-\eta_0-\frac52}T^{\frac52+\epsilon}\Big)\notag\\
&\ll_\epsilon y^{\frac78}T^{-\frac12+\frac{11}2\eta+4\epsilon}+y^{\frac74}T^{-1+\frac{11}2\eta+4\epsilon}.
\end{align}
When counting the quantity $\ll q^4\log^4T$, we observe that these $Z_1$ contribute an error
\begin{align}
\ll_\epsilon y^{\frac78}T^{-\frac12+\frac{11}2\eta+\epsilon}+y^{\frac74}T^{-1+\frac{11}2\eta+\epsilon}
\end{align}
to $Z$. Then (B1) of the proposition follows from this and the estimate on the contribution of all $Z_1$ in case one. Thus we complete our proofs of the proposition and Theorem \ref{thm1}.

\section{Application to the proportion of critical zeros}\label{sec5}
Suppose that $\log q=o(\log T)$ and $\chi$ is a primitive Dirichlet character with $q$ its modulus. Let $Q(x)$ be a real polynomial satisfying $Q(0)=1$ and $Q'(x)=Q'(1-x)$.
The well-known Levinson method for the Dirichlet $L$-function yields the inequality (see also \cite{Con-More} and Appendix A in \cite{CIS-Crit})
\begin{align}\label{+5.2}
\kappa(\chi)\ge1-\frac1R\log(T^{-1}I_R(Q,\chi))+o(1)
\end{align}
for any given positive real constant $R$. Actually, if $Q(x)$ is a linear polynomial, the inequality gives a lower bound for the proportion of simple zeros $\kappa^*(\chi)$.

To specify the mollifier in \eqref{xm1.18}, we take
\begin{align}
B(s,\chi)=\sum_{n\le y}\frac{\chi(n)a(n)}{n^{s}},
\end{align}
where $y=T^\theta$ with $\theta=\frac47-\epsilon$ and
\begin{align}\label{5.2}
a(n)=\mu(n)\bigg(P_1\bigg(\frac{\log y/n}{\log y}\bigg)+P_2\bigg(\frac{\log y/n}{\log y}\bigg)\sum_{p\mid n,~p\le y^{3/4}}P\bigg(\frac{\log p}{\log y}\bigg)\bigg)
\end{align}
with $P_1,~P_2,~P$ are real polynomials that meet $P_1(0)=P_2(0)=P(0)=0$ and $P_1(1)=1$.

Since $\alpha,~\beta\ll\mathcal{L}^{-1}_\chi$, $\big(\frac{2\pi}{qT}\big)^{\alpha+\beta}=e^{-a-b}$ and
\begin{align}
L(1+s,\chi_0)=\frac{\phi(q)}{q}\bigg(\frac1s+\gamma+c_q+o(1)\bigg)
\end{align}
as $s\rightarrow0$, it follows directly from \eqref{+1.10} that
\begin{align}\label{x6.5}
I_R(Q,\chi)\sim T\frac{\phi(q)}{q}Q\bigg(\frac{-d}{da}\bigg)\overline{Q}\bigg(\frac{-d}{db}\bigg) \frac{\Sigma(\beta,\alpha)-e^{-a-b} \Sigma(-\alpha,-\beta)}{\alpha+\beta}\bigg|_{a=b=-R}
\end{align}
with
\begin{align}\label{+1.11}
\Sigma(\alpha,\beta)=\sum_{h,k\le y}\frac{a(h)a(k)}{h^{1+\alpha}k^{1+\beta}}\chi_0(hk)(h,k)^{1+\alpha+\beta}.
\end{align}
Let
\begin{align}\label{5.3}
F(j,s)=\prod_{p\mid j}\bigg(1-\frac{1}{p^s}\bigg),
\end{align}
then we have
\begin{align}\label{5.4}
\Sigma(\alpha,\beta)=\sum_{j\le y}j^{-1}F(j,1+\alpha+\beta)\chi_0(j)E(\alpha,j)E(\beta,j)
\end{align}
with
\begin{align}\label{5.5}
E(\alpha,j)=\sum_{h\le y/j}\frac{a(hj)\chi_0(h)}{h^{1+\alpha}}.
\end{align}

To estimate $\Sigma(\alpha,\beta)$, we need following lemmas, which may be proved similarly as corresponding results for the Riemann zeta-function, see also \cite{Bau,Con-Zero, Lev-More,Sel}.
\begin{lemma}\label{lemma5.1}
Let $\chi_0$ be the principle character mod $q$ and $P$ be a real polynomial with $P(0)=0$.  Suppose that
\begin{align}
G=\sum_{{n\le y/j}\atop{(n,j)=1}}\frac{\mu(n)}{n^{1+\alpha}}\chi_0(n)P\bigg(\frac{\log y/nj}{\log y}\bigg).
\end{align}
Then we have
\begin{align}
G=&\frac1{F(qj,1+\alpha,\chi_0)}\bigg(\alpha P\bigg(\frac{\log y/j}{\log y}\bigg)+\frac1{\log y}P'\bigg(\frac{\log y/j}{\log y}\bigg)\bigg)\notag\\
&+O\bigg(\frac{(\log\log y)^3F_1(qj,1-2\delta)}{\log^2 y}\bigg)+O\bigg(\frac{(\log\log y)^3F_1(qj,1-2\delta)}{\log y}\bigg(\frac jy\bigg)^{\delta/M}\bigg)
\end{align}
uniformly in $j\le y,~\log q=o(\log y)$ and $\alpha\ll\frac1{\log y}$. Here $F(j,s)$ is defined by \eqref{5.3}, $F_1(j,s)=\prod_{p\mid d}(1+1/p^{-s})$, $\delta=1/\log\log y$ and $M$ is an absolute constant.
\end{lemma}

\begin{proof}
Since $\chi_0(n)=1$ for $(n,1)=1$ and vanishes otherwise, we note that
\begin{align}
G=\sum_{{n\le y/j}\atop{(n,qj)=1}}\frac{\mu(n)}{n^{1+\alpha}}P\bigg(\frac{\log y/nj}{\log y}\bigg).
\end{align}
Then one may prove the lemma the same as Lemma 10 in \cite{Con-Zero} with $F(j,s)$ replaced by $F(qj,s)$.
\end{proof}

\begin{lemma}\label{lemma5.2}
Let $f(p)=1+O(p^{-c}),~c>0$ and $f(n)=\prod_{p\mid n}f(p)$. If $P$ is a polynomial and $\chi_0$ denotes the principal character modulo $q\ge1$, we have, for any integer $k\ge0$,
\begin{align}
\sum_{n\le y}\frac{\mu^2(n)\chi_0(n)}{n}f(n)\log^k\frac yn=\frac{\phi(q)}{q}\cdot \mathcal{P}_f \cdot(k+1)^{-1}\log^{k+1}y+O_k(\log^k y),
\end{align}
where
\begin{align}
\mathcal{P}_f=\prod\bigg(1+\frac{\chi_0(p)(f(p)-1)}{p+1}\bigg)\bigg(1-\frac{\chi_0(p)}{p^2}\bigg).
\end{align}
In the special case
\begin{align}
f(p)=\frac{1-p^{-1-\alpha-\beta}}{(1-p^{-1-\alpha})(1-p^{-1-\beta})},
\end{align}
it follows that $\mathcal{P}_f=1+O(\mathcal{L}^{-1}_\chi)$ for $\alpha, \beta\ll\mathcal{L}^{-1}_\chi$.
\end{lemma}
\begin{proof}
Regarding the term $\mu^2(n)\chi_0(n)$ we follow the procedure of the proof of Lemma 3.11 in \cite{Lev-More},  which deals with $\mu^2(n)$. Then we prove the $k=0$ case almost the same as in \cite{Lev-More}. The only difference is due to the following equation
\begin{align}
\sum_{r\le y}\frac{\chi_0(r)}r=\frac{\phi(q)}{q}\log y+O(1)
\end{align}
used at the end. For $k\ge1$ it follows by using Abel summation with the $k=0$ case.
\end{proof}

\begin{lemma}\label{lemma5.3}
Let $\delta,~\delta'\ge0$ and $\delta+\delta'\le c<1$. Further, let $F_1(j,s)=\prod_{p\mid d}(1+1/p^{-s})$ as before. Then, for any positive integer $r$, we have
\begin{equation}
\sum_{{n\le y}\atop{(n,q)=1}}\frac{\mu^2(n)}{n^{1-\delta}}F_1(n,1-\delta')^r=\left\{
\begin{aligned}
&O_{c,r}(\log y)\ \ \ &\text{if}~\delta=0,\\
&O_{c,r}(y^\delta/\delta)\ \ \ &\text{if}~\delta>0.
\end{aligned}
\right.
\end{equation}
\end{lemma}
\begin{proof}
For $\delta=0$, it is a direct result of Lemma \ref{lemma5.2}; For $\delta>0$, it is a trivial bound.
\end{proof}

We also need following two well-known results.
\begin{lemma} [Mertens Theorem]
\label{Mertems}
\begin{align}
\sum_{p\le y}\frac{\log p}p=\log y+O(1).
\end{align}
\end{lemma}

\begin{lemma} [Levinson \cite{Lev-More}]
\label{Levinson}
\begin{align}
\sum_{p\mid n}\frac{\log p}p\ll\log\log n.
\end{align}
\end{lemma}
Two estimates
\begin{align}
q^{1+\alpha}=q(1+o(1))\ \ \ \ \text{and} \ \ \ \ \prod_{p\mid q}\bigg(1-\frac1{p^{1+\alpha}}\bigg)=\frac{\phi(q)}{q}(1+o(1))
\end{align}
as $\alpha\rightarrow0$ or $T\rightarrow\infty$ will be used frequently in our following calculation, and we will not point them out especially. The first estimate is obvious, and the second one follows from Lemma \ref{Levinson} by taking the logarithm.

In the expression of $E(\alpha,j)$ with $a(hj)$ given by \eqref{5.2}, we separate the sum $\sum_{p\mid hj}$ to $\sum_{p\mid j}+\sum_{p\mid h}$ and have
\begin{align}
E(\alpha,j)=&\mu(j)\sum_{{h\le y/j}\atop{(h,j)=1}}\frac{u(h)\chi_0(h)}{h^{1+\alpha}}\bigg(P_1\bigg(\frac{\log y/hj}{\log y}\bigg)+P_2\bigg(\frac{\log y/hj}{\log y}\bigg)\sum_{p\mid hj,~p\le y^{3/4}}P\bigg(\frac{\log p}{\log y}\bigg)\bigg)\notag\\
=&\mu(j)\sum_{{h\le y/j}\atop{(h,j)=1}}\frac{u(h)\chi_0(h)}{h^{1+\alpha}}P_1\bigg(\frac{\log y/hj}{\log y}\bigg)\notag\\
&+\mu(j)\sum_{{p\le \min(y/j,~y^{3/4})}\atop{(p,qj)=1}}\frac{\mu(p)P\Big(\frac{\log p}{\log y}\Big)}{p^{1+\alpha}}\sum_{{h\le y/pj}\atop{(h,pj)=1}}\frac{u(h)\chi_0(h)}{h^{1+\alpha}}P_2\bigg(\frac{\log y/hj}{\log y}\bigg)\notag\\
&+\mu(j)\sum_{{p\mid i}\atop{p\le y^{3/4}}}P\bigg(\frac{\log p}{\log y}\bigg)\sum_{{h\le y/j}\atop{(h,j)=1}}\frac{u(h)\chi_0(h)}{h^{1+\alpha}}P_2\bigg(\frac{\log y/hj}{\log y}\bigg).
\end{align}
Thus we denote
\begin{align}
E(\alpha,j)=E_1(\alpha,j)+E_2(\alpha,j)+E_3(\alpha,j)
\end{align}
with obvious meanings. Let
\begin{align}
\mathcal{V}_i(a,t)=a\theta P_i(t)+p'_i(t)
\end{align}
with $i=1,~2$ and
\begin{align}
\mathcal{W}_2(a,t)=\int_0^{\min(t,\frac34)} e^{-a\theta u}\frac{P(u)}{u}\mathcal{V}_2(a,t-u)du.
\end{align}
Using Lemma \ref{lemma5.1} we get
\begin{align}\label{+5.21}
E_1(\alpha,j)=&\mathcal{L}^{-1}\frac{\mu(j)}{F(qj,1+\alpha)}\mathcal{V}_1\bigg(a,\frac{\log y/j}{\log y}\bigg)+O(A_1)+O(B_1),
\end{align}
\begin{align}\label{+5.22}
E_2(\alpha,j)=&-\mathcal{L}^{-1}\frac{\mu(j)}{F(qj,1+\alpha)}\mathcal{W}_2\bigg(a,\frac{\log y/j}{\log y}\bigg)+O(A_2)+O(B_2)
\end{align}
and
\begin{align}\label{+5.23}
E_3(\alpha,j)=&\mathcal{L}^{-1}\frac{\mu(j)}{F(qj,1+\alpha)}\mathcal{V}_2\bigg(a,\frac{\log y/j}{\log y}\bigg)\sum_{{p\mid j}\atop{p\le y^{3/4}}}P\bigg(\frac{\log p}{\log y}\bigg)+O(A_3)+O(B_3)
\end{align}
with
\begin{align}
A_i\ll\frac{(\log\log y)^4F_1(qj,1-2\delta)}{\log^2 y}\Bigg(1+\sum_{p\mid j}P\bigg(\frac{\log p}{\log y}\bigg)+\sum_{p\le y/j}\frac{F_1(p,1-2\delta)P\Big(\frac{\log p}{\log y}\Big)}{p}\Bigg)
\end{align}
and
\begin{align}
B_i\ll\frac{(\log\log y)^4F_1(qj,1-2\delta)}{\log y}\bigg(\frac jy\bigg)^{\delta/M}\Bigg(1+\sum_{p\mid j}P\bigg(\frac{\log p}{\log y}\bigg)+\sum_{p\le y/j}\frac{F_1(p,1-2\delta)P\Big(\frac{\log p}{\log y}\Big)}{p^{1-\delta/M}}\Bigg)
\end{align}
for all $i=1,2$ and $3$ since $1\ll p^\alpha\ll1$ and $F_1(pj,1-2\delta)=F_1(p,1-2\delta)F_1(j,1-2\delta)$ for $(p,j)=1$. In the evaluation of $E_2$ above, we firstly use $\sum_{(p,qj)=1}=\sum_p-\sum_{p\mid qj}$ to remove the co-primality condition for the sum on $p$, then we estimate the first term  by Lemma \ref{Mertems} and move the second term to error terms due to the fact that $\sum_{p\mid qj}\log p/p\ll\log\log y$.

We employ \eqref{+5.21}-\eqref{+5.23} into \eqref{5.4} to separate $\Sigma(\alpha,\beta)$ to some terms with obvious meanings. Each term contains a sum on the variable $j$ and at most two prime variables. We firstly consider the terms containing $A_i$ or $B_i$. If a term does not contain any prime variable in its sum, we can estimate it directly by Lemma \ref{lemma5.3}; Otherwise, by interchanging the order we can make it true that the innermost sum is on $j$ and the other is on prime variables. As $j$ is square free, if a term contains two prime variables in its sum, we should employ the formula
\begin{align}
\sum_{p_1\mid j}\sum_{p_2\mid j}=\sum_{p\mid j}+\sum_{{p_1p_2\mid j}\atop{p_1\neq p_2}}
\end{align}
before interchanging the order.  Since $(q,j)=1$, we have
\begin{align}
F(qj,1+\alpha)=F(q,1+\alpha)F(j,1+\alpha)  \ \ \ \ \text{and}  \ \ \ \ \ F(q,1+\alpha)=\frac{\phi(q)}{q}(1+\mathcal{L}^{-1}_\chi).\notag
\end{align}
As $F(j,1+\alpha)^{\pm1}=O(F_1(j,1-\delta))$ for all $\alpha\ll1/\mathcal{L}_\chi$, using Lemma \ref{lemma5.3} in the innermost sum and employing
\begin{align}
\sum_{p\le y}\frac{P\Big(\frac{\log p}{\log y}\Big)}{p}\ll1 \ \ \ \text{and}~ \sum_{p\le y}\frac{F_1(p,1-2\delta)P\Big(\frac{\log p}{\log y}\Big)}{p}\ll1,
\end{align}
which can be deduced from Lemma \ref{Mertems}, to the sum on prime variables, we find that these terms contribute an error $\ll(\log\log y)^7\log^{-2}y$ to $\Sigma(\alpha,\beta)$.

By employing main terms of $E_i(\alpha,j)$ and $E_i(\beta,j)$ in the sum of $\Sigma(\alpha,\beta)$ in \eqref{5.4} we have
\begin{align}
\Sigma(\alpha,\beta)\sim\frac{q^2}{\phi(q)^2}\mathcal{L}^{-2}\sum_{j\le y}\frac{\mu^2(j)\chi_0(j)F(j,1+\alpha+\beta)}{jF(j,1+\alpha)F(j,1+\beta)}\mathcal{G}(\alpha,y,j) \mathcal{G}(\beta,y,j),
\end{align}
where
\begin{align}
\mathcal{G}(\alpha,y,j)=\mathcal{V}_1\bigg(a,\frac{\log y/j}{\log y}\bigg)+\mathcal{V}_2\bigg(a,\frac{\log y/j}{\log y}\bigg)\sum_{{p\mid j}\atop{p\le y^{3/4}}}P\bigg(\frac{\log p}{\log y}\bigg)+\mathcal{W}_2\bigg(a,\frac{\log y/j}{\log y}\bigg),
\end{align}
and $\mathcal{G}(\beta,y,j)$ has a similar expression. Then we can split $\Sigma(\alpha,\beta)$ into nine terms and evaluate them term by term. If a term does not contain any prime variable in its sum, we evaluate it directly by Lemma \ref{lemma5.2}; If a term contains prime variables in its sum, by interchanging the order we can also make it true that the innermost sum is on $j$ and the other is on prime variables. Also, the formula
\begin{align}
\sum_{p_1\mid j}\sum_{p_2\mid j}=\sum_{p\mid j}+\sum_{{p_1p_2\mid j}\atop{p_1\neq p_2}}
\end{align}
should be employed first when a term contains two prime variables in its sum. Then we evaluate the innermost sum on $j$ by Lemma \ref{lemma5.2} with $q$ replaced by $qp$ or $qp_1p_2$ according to the number of prime variables it contains. The outer sum can be estimated by Lemma \ref{Mertems} according to prime variables one by one. After doing these we actually have
\begin{align}\label{5.31}
\Sigma(\alpha,\beta)=\sum_{1\le i_1,i_2\le3}\Sigma_{i_1i_2}(a,b)+O\Big((\log\log y)^7\log^{-2}y\Big),
\end{align}
where $\Sigma_{i_1i_2}(a,b)$ denotes the main term of the sum on $E_{i_1}(\alpha,j)E_{i_2}(\beta,j)$, and having that $\Sigma_{i_1i_2}(a,b)=\Sigma_{i_2i_1}(b,a)$ with
\begin{align}\label{+5.31}
\Sigma_{11}(a,b)=\frac{q}{\phi(q)}\cdot\frac1{\theta\mathcal{L}}\int_0^1\mathcal{V}_1(a,1-t) \mathcal{V}_1(b,1-t)dt,
\end{align}

\begin{align}\label{+5.32}
\Sigma_{12}(a,b)=-\frac{q}{\phi(q)}\cdot\frac1{\theta\mathcal{L}}\int_0^1\mathcal{V}_1(a,1-t) \mathcal{W}_2(b,1-t)dt,
\end{align}

\begin{align}\label{+5.33}
\Sigma_{22}(a,b)=-\frac{q}{\phi(q)}\cdot\frac1{\theta\mathcal{L}}\int_0^1\mathcal{W}_2(a,1-t) \mathcal{W}_2(b,1-t)dt,
\end{align}

\begin{align}\label{+5.34}
\Sigma_{13}(a,b)=\frac{q}{\phi(q)}\cdot\frac1{\theta\mathcal{L}}\int_0^{\frac34}\frac{P(t)}{t}dt \int_0^{1-t}\mathcal{V}_1(a,1-t-t_1) \mathcal{V}_2(b,1-t-t_1)dt_1,
\end{align}

\begin{align}\label{+5.35}
\Sigma_{23}(a,b)=-\frac{q}{\phi(q)}\cdot\frac1{\theta\mathcal{L}}\int_0^{\frac34}\frac{P(t)}{t}dt \int_0^{1-t}\mathcal{W}_2(a,1-t-t_1) \mathcal{V}_2(b,1-t-t_1)dt_1
\end{align}
and
\begin{align}\label{+5.36}
&\Sigma_{33}(a,b)=\frac{q}{\phi(q)}\cdot\frac1{\theta\mathcal{L}}\Bigg\{\int_0^{\frac34}\frac{P(t)^2}{t}dt \int_0^{1-t}\mathcal{V}_2(a,1-t-t_1) \mathcal{V}_2(b,1-t-t_1)dt_1\notag\\
&+\int_0^{\frac34}\frac{P(t)}{t}dt\int_0^{\min(\frac34,1-t)}\frac{P(t_1)}{t_1}dt_1 \int_0^{1-t-t_1}\mathcal{V}_2(a,1-t-t_1-t_2) \mathcal{V}_2(b,1-t-t_1-t_2)dt_2\Bigg\}.
\end{align}

\subsection{Proof of Theorem \ref{thm2}} We have from \eqref{+5.2} and \eqref{x6.5} that
\begin{align}
\kappa\ge1-\frac1R\log \bigg(Q\bigg(\frac{-d}{da}\bigg)Q\bigg(\frac{-d}{db}\bigg) \frac{\phi(q)}q\frac{\Sigma(\beta,\alpha)-e^{-a-b}\Sigma(-\alpha,-\beta)} {\theta(a+b)}\bigg)\bigg|_{a=b=-R}.
\end{align}
Then we substitute main terms of \eqref{+5.31}-\eqref{+5.36} into the above formula by \eqref{5.31}, also use Mathematica with the following choices of parameters. With $\theta=\frac47-\epsilon$, $R=1.3$,
\begin{align}
&Q(x)=1-.642x-1.227(x^2/2-x^3/3)-5.178(x^3/3-x^4/2+x^5/5),\notag\\
&P_1(x)=x-.617x(1-x)-.125x^2(1-x)-.148x^3(1-x),\notag\\
&P_2(x)=x,\notag\\
&P(x)=1.155x-1.564x^2+.177x^3,\notag
\end{align}
we have $\kappa(\chi)>.4172$. To get $\kappa^*(\chi)>.4074$ we take $R=1.116$,
\begin{align}
&Q(x)=1-1.032x,\notag\\
&P_1(x)=x-.525x(1-x)-.183x^2(1-x)-.085x^3(1-x), \ \ \ \ \ \ \ \ \ \ \ \ \ \ \ \ \ \ \notag\\
&P_2(x)=x,\notag\\
&P(x)=.838x-.938x^2-.084x^3.\notag
\end{align}
Thus we prove Theorem \ref{thm2}.\\

\noindent
{\bf Remark.}
If we take the coefficient in a more general form
\begin{align}\label{x6.43}
a(n)=\mu(n)\bigg(P_0\bigg(\frac{\log y/n}{\log y}\bigg)+\sum_{1\le k\le K}P_k\bigg(\frac{\log y/n}{\log y}\bigg)\sum_{p\mid n,~p\le y^{3/4}}\bigg(\frac{\log p}{\log y}\bigg)^k\bigg)
\end{align}
with $K=3$, we can improve numeric results slightly to
\begin{align}
\kappa(\chi)>.417277\ \ \ \ \ and \ \ \ \ \ \kappa^*(\chi)>.407475.
\end{align}
Moreover, if removing the condition $p\le y^{\frac34}$ in the sum of \eqref{x6.43} with the help of \cite{PRZZ}, we may have
\begin{align}
\kappa(\chi)>.417293\ \ \ \ \ and \ \ \ \ \ \kappa^*(\chi)>.40751.
\end{align}

\section*{Acknowledgements}
\noindent
We would like to express our heartfelt thanks to the anonymous referee for his careful reading and helpful suggestion.


\begin{thebibliography}{99}

\bibitem{BCH} R. Balasubramanian, J. B. Conrey and D. R. Heath-Brown, Asymptotic mean square of the product of the Riemann zeta-function and a Dirichlet polynomial, J. Reine Angew. Math., 357 (1985), 161-181.
\bibitem{Bau} P. Bauer, Zeros of Dirichlet $L$-series on the critical line, Acta Arith., 93 (2000), 37-52.
\bibitem{BG} S. Bettin and S. Gonek, The $\theta=\infty$ conjecture implies the Riemann hypothesis, Mathematika, 63 (2017) no. 1, 29-33.
\bibitem{BC} S. Bettin and V. Chandee, Trilinear forms with kloosterman fractions, Adv. Math., 328 (2018), 1234-1262.
\bibitem{BCR} S. Bettin, V. Chandee and M. Radziwi{\l}{\l}, The mean square of the product of the Riemann zeta-function with Dirichlet polynomials, J. Reine Angew. Math., 720 (2017), 51-79.
\bibitem{BCY} H. Bui, J. B. Conrey and M. Young, More than 41\% of the zeros of the zeta function are on the critical line, Acta Arith., 150 (2011), no.1, 35-64.
\bibitem{Bui} H. Bui, Critical zeros of the Riemann zeta-function, available at arXiv: 1410.2433.
\bibitem{Con-Zero} J. B. Conrey, Zeros of derivatives of the Riemann's $\xi$-function on the critical line, J. Number Theory, 16 (1983), 49-74.
\bibitem{CG} J. B. Conrey and A. Ghosh, A simpler proof of Levinson's theorem, Math. Proc. Cambridge Philos. Soc., 97 (1985), 385-395.
\bibitem{CGG} J. B. Conrey, A. Ghosh and S. M. Gonek, Large gaps between zeros of the zeta-function, Mathematika, 33 (1986), no. 2, 212-238.
\bibitem{Con-More} J. B. Conrey, More than two fifths of the zeros of the Riemann zeta function are on the critical line, J. Reine Angew. Math., 339 (1989), 1-26.
\bibitem{CIS-Crit} J. B. Conrey, H. Iwaniec and K. Soundararajan, Critical zeros of Dirichlet $L$-function, J. Reine Angew. Math., 681 (2013), 175-198.
\bibitem{CIS-Asym} J. B. Conrey, H. Iwaniec and K. Soundararajan, Asymptotic large sieve, available at arXiv: 1105.1176.
\bibitem{DIK} J. -M. Deshouillers and H. Iwaniec, Kloosterman sums and Fourier coefficients of cusp forms, Invent. Math., 70 (1982), 219-288.
\bibitem{DIP} J. -M. Deshouillers and H. Iwaniec, Power mean values of the Riemann zeta function II, Acta Arith., 48 (1984), 305-312.
\bibitem{Feng} S. Feng, Zeros of the Riemann zeta function on the critical line, J. Number Theory, 132 (4) (2012), 511-542.
\bibitem{Hea} D. R. Heath-Brown, Simple zeros of the Riemann zeta-function on the critical line, Bull. Lond. Math. Soc., 11 (1979), 17-18.
\bibitem{HeaV} D. R. Heath-Brown, Prime numbers in short intervals and a generalized Vaughan identity, Canad. J. Math., 34 (1982), no. 6, 1365-1377.
\bibitem{KRZ} P. K\"{u}hn, N. Robles and D. Zeindler, On mean values of mollifiers and $L$-functions associated to primitive cusp form, Math. Z., (2018), available at https://doi.org/10.1007/s00209-018-2099-9.
\bibitem{Lev-More} N. Levinson, More than one third of zeros of Riemann's zeta-function are on $\sigma=1/2$, Adv. Math., 13 (1974), 383-436.
\bibitem{LR} X. Li and M. Radziwi{\l}{\l}, The Riemann-zeta function on vertical arithmetic progressions, Int. Math. Res. Not. IMRN, (2015), no. 2, 325-354.
\bibitem{PR} K. Pratt and N. Robles, Perturbed moments and a longer mollifier for critical zeros of $\zeta$, Res. Number Theory 4(2018), no. 1, Art. 9, 26pp.
\bibitem{Rad} M. Radziwi{\l}{\l}, Limitations to mollifying $\zeta(s)$, available at arXiv: 1207.6583.
\bibitem{RRZ} N. Robles, A. Roy and A. Zaharescu, Twisted second moments of the Riemann zeta-function and applications, J. Math. Anal. Appl., 434 (2016), no. 1, 271-314.
\bibitem{PRZZ} K. Pratt, N. Robles, A. Zaharescu and D. Zeindler, Conbinatorial applications of autocorrelation ratios, available at arXiv: 1802.10521.
\bibitem{Sel} A. Selberg, On the zeros of Riemann's zeta-function, Skr. Norske Videnskaps-Akad. Oslo, I, 10 (1942), 1-59.
\bibitem{Son} K. Sono, An application of generalized mollifiers to the Riemann zeta-function, Kyushu J. Math., 72 (2018), 35-69.
\bibitem{Sou} K. Soundararajan, Mean-values of the Riemann zeta-function, Mathematika, 42 (1995), no. 1, 158-174.
\bibitem{Wu} X. Wu, Distinct zeros and simple zeros for the family of Dirichlet $L$-functions, Quart. J. Math., 67 (2016), 757-779.
\bibitem{You} M. P. Young, A short proof of Levinson's theorem, Arch. Math., 95 (2010), 539-548.










\end{thebibliography}
\end{document}